\documentclass[letterpaper]{article}

\usepackage[margin=2.5cm]{geometry}
\usepackage{amsmath,amsthm,amsfonts,amssymb, tikz, array, tabularx}
\usetikzlibrary{calc}

\newcommand{\overbar}[1]{\mkern 1.5mu\overline{\mkern-1.5mu#1\mkern-1.5mu}\mkern 1.5mu}

\newcommand{\ex}{\mbox{ex}}

\newtheorem{thm}{Theorem}[section]
\newtheorem{lemma}[thm]{Lemma}
\newtheorem{prop}[thm]{Proposition}

\theoremstyle{definition}

\newtheorem{clm}[thm]{Claim}
\newtheorem{defn}[thm]{Definition}

\numberwithin{equation}{section}
\numberwithin{figure}{section}

\title{Triangles in graphs without bipartite suspensions\footnotetext[0]{Accepted version for publication at \emph{Discrete Mathematics}.}}

\author{
	Dhruv Mubayi\footnote{Department of Mathematics, Statistics and Computer Science, University of Illinois, Chicago, IL 60607, USA. \texttt{Email:mubayi@uic.edu}.} $\qquad$ 
	Sayan Mukherjee\footnote{Blueqat Research, Blueqat Inc., Tokyo 150-6139, Japan.  \texttt{Email:sayan@blueqat.com}} \footnote{Department of Physics, The University of Tokyo, Tokyo 113-0033, Japan.}
}

\begin{document}
	
	\maketitle
	
	\begin{abstract}
		Given graphs $T$ and $H$, the generalized Tur\'an number $\ex(n,T,H)$ is the maximum number of copies of $T$ in an $n$-vertex graph with no copies of $H$. Alon and Shikhelman, using a result of Erd\H os, determined the asymptotics of $\ex(n,K_3,H)$ when the chromatic number of $H$ is greater than three and proved several results when $H$ is bipartite. We consider this problem when $H$ has chromatic number three. Even this special case for the following relatively simple three chromatic graphs appears to be challenging. The suspension $\widehat H$ of a graph $H$ is the graph obtained from $H$ by adding a new vertex adjacent to all vertices of $H$. We give new upper and lower bounds on $\ex(n,K_3,\widehat{H})$ when $H$ is a path, even cycle, or complete bipartite graph. One of the main tools we use is the triangle removal lemma, but it is unclear if much stronger statements can be proved without using the removal lemma.
	\end{abstract}

	
	\section{Introduction}
	
	A graph $G=(V(G),E(G))$ consists of a vertex set $V(G)$ and an edge set $E(G)\subseteq \binom{V(G)}{2}$. Let $e(G)$ denote the number of edges of $G$. Say that $G$ is $F$-free if $G$ contains no subgraph isomorphic to $H$. We emphasize that we do not consider induced subgraphs in this definition.
	
	For graphs $T$ and $H$ with no isolated vertices and integer $n$, the generalized Tur\'an number $\ex(n,T,H)$ is the largest number of copies of $T$ in an $H$-free $n$-vertex graph. When $T=K_2$, this is the  Tur\'an number $\ex(n,H)$ of the graph $H$.
	
	The systematic study of $\ex(n, T, H)$ for $T\neq K_2$ was initiated by Alon and Shikhelman~\cite{alonTCopies2016}. Before then, there had been sporadic results determining this function for several $T$ and $H$, beginning with $\ex(n,K_t,K_r)$ for $t<r$ (see \cite{bollobasKr1976,erdosKtKr1962}). Several cases where $H=K_r$ were studied in \cite{gyoriKrfree1991}. There has been a lot of recent activity when $(T,H)=(K_3,C_{2k+1})$ and $(T,H)=(C_5,K_3)$ \cite{alonTCopies2016, bollobasGyoriK3C52008, methukuK3C52019, gyoriK3OddCycle2012, grzesikC5K32012}. In \cite{gerbnerGeneral2019}, the cases $(T,H)=(P_k,K_{2,t})$ and $(T,H)=(C_k,K_{2,t})$ have also been studied and some generic bounds on $\ex(n,T,H)$ are given. See also \cite{luoKrLongCycle2018} for a related result about the number of $s$-cliques in graphs without cycles of length at least $k$.
		
	Alon and Shikhelman \cite{alonTCopies2016} determine all pairs of graphs $T,H$ with $\ex(n,T,H)=\Theta(n^{|V(T)|})$. Further, they prove that if $T$ and $H$ are trees, then there exists an $m(T,H)$ such that $\ex(n,T,H)=\Theta(n^{m(T,H)})$.
	They also study the problem when $H$ is a tree and $T$ is a bipartite graph, and give several results on $\ex(n,K_t,H)$ for bipartite $H$. One general result they prove using a theorem of Erd\H os \cite{erdosComplete-t-partite1964} is that if the chromatic number $\chi(H)>t$, then
	\[
	\ex(n,K_t,H)=\binom{\chi(H)-1}t \left(\frac n{\chi(H)-1}\right)^t+o(n^t).
	\]
 In \cite{jieSharpResults2020}, the error term was determined more precisely.
	
All our results concern $T=K_3$. Since the asymptotic formula of $\ex(n,K_3,H)$ is already known for $\chi(H)>3$ and \cite{alonTCopies2016} studies the case $\chi(H)=2$ quite extensively, we consider the wide open case $\chi(H)=3$. Even within this class, we restrict our attention to a very specific and simple family of $3$-chromatic graphs. For any graph $H$, let $\widehat{H}$ denote the suspension $K_1\vee H$ obtained by adding a new vertex adjacent to every vertex of $H$. We  obtain upper and lower bounds on  $\ex(n,K_3,\widehat{H})$ for
$H \in \{K_{a,b}, C_{2k}, P_k\}$, where $P_k$ denotes the path on $k$ edges.
	
	Given a graph $G=(V,E)$ and a vertex $v\in V$, let $N_G(v)=\{u\in V: uv\in E\}$ denote the neighborhood of $v$ in $G$. For any subset $X\subseteq V$, let $e(X)$ denote the number of edges in the subgraph $G[X]$ induced by $X$.  Let $t(G)$ denote the number of triangles in $G$.
	
	If $G$ is $\widehat H$-free, then $N_G(v)$ is  $H$-free for every $v \in V$. This implies that 
	\[
	t(G)=\frac 13\sum_{v\in V}e(N_G(v))\le \frac13 \sum_{v\in V}\ex(|N_G(v)|,H)\le \frac n3\cdot \ex(n,H).
	\]
	Hence,
	\begin{equation} 
	\label{eq:introupper}
	\ex(n,K_3,\widehat H)\le \frac n3\cdot \ex(n,H).
	\end{equation}
	
	\medskip
	All our results give improvements on (\ref{eq:introupper}). For our first result $H=K_{a,b}$ and $\widehat H = K_{1,a,b}$, where $1\le a\le b$. Here (\ref{eq:introupper}) combined with the K\"{o}vari-S\'{o}s-Tur\'{a}n theorem \cite{kovariSosTuran1954}, which asserts that $\ex(n,K_{a,b})=O(n^{2-\frac 1a})$ yields $\ex(n, K_3, K_{1,a,b}) = O(n^{3-\frac1a})$. We improve this as follows.
	\begin{thm}
		\label{thm:k1ab}
		For fixed $1\le a\le b$ and $n\to\infty$,
		\begin{equation}
		\label{eq:k1ab}
		\mbox{\rm \ex}(n,K_3,K_{1,a,b}) = o(n^{3-\frac 1a}).
		\end{equation}
	\end{thm}
	
	Notice that setting $a=b=2$ in (\ref{eq:k1ab}) yields $\ex(n,K_3,K_{1,2,2})=o(n^{5/2})$, where $K_{1,2,2} = 
	\widehat{C}_4$. This is related to a question in \cite{mubayiVerstraeteSurveyExpansion2016}, where the authors asked whether $\ex(n,K_3,K_{1,2,2})=O(n^2)$. While this remains open we do give a quadratic lower bound  in Proposition \ref{prop:w4Lowerbd}. Narrowing the (huge) gap in the bounds $$\Omega(n^2)<\ex(n, K_3, K_{1,2,2}) < o(n^{5/2})$$
	 is perhaps the most basic and attractive open problem raised in this paper.
	
	\medskip
	Our next result concerns  $H=C_{2k}$. Here  (\ref{eq:introupper}) together with  the classical bound $\ex(n,C_{2k})=O(n^{1+\frac 1k})$ of Bondy-Simonovits \cite{bondySimonovits1974}  yields $\ex(n,K_3,\widehat{C}_{2k}) = O(n^{2+\frac 1k})$.
	\begin{thm}
		\label{thm:w2k}
		For fixed $k\ge 2$ and $n\to\infty$,
		\begin{equation}
		\label{eq:w2k}
		\ex(n,K_3,\widehat{C}_{2k}) = o(n^{2+\frac 1k}).
		\end{equation}
	\end{thm}
	
	The same lower bound construction for $\ex(n, K_3, K_{1,2,2})$ yields $\ex(n, K_3, \widehat{C}_{2k}) =\Omega( n^2)$.
	
	 Our final results concern $\ex(n,K_3,\widehat{P}_k)$ for $k\ge 3$.
	 We begin with a simple proposition. 
	
	\begin{prop}
		\label{thm:pkhat}
		Let $n \ge k\ge 3$. Then 
		\begin{equation}
		\label{eq:thm-pkhat-upper}
		\left\lfloor\frac{k-1}{2}\right\rfloor\cdot \frac{n^2}{8}\le \ex(n,K_3,\widehat{P}_k)\le  \frac{k-1}{12}\cdot n^2+\frac{(k-1)^2}{12}\cdot n,
		\end{equation}
		where the lower bound holds when $n$ is a multiple of $4\lfloor \frac{k-1}{2}\rfloor$.
	\end{prop}
	We believe that the lower bound above is asymptotically tight for all fixed $k \ge 3$ and prove this for the first three cases $k=3$, $4$ and $5$.
	
	\begin{thm}
		\label{thm:p3p4p5hat}
		For  $k=3$, $4$ and $5$, 
		\begin{equation}
		\label{eq:thm-p3p4hat}
		\ex(n,K_3,\widehat{P}_k)= \left\lfloor\frac{k-1}{2}\right\rfloor\cdot \frac{n^2}{8}+o(n^2).
		\end{equation}
				When $k=3$ or $k=5$, the error term can be improved to $O(n)$. 	
	\end{thm}
	
	In Section 2, we present some preliminary results that we use in our proofs. In Section 3 we  prove Theorem \ref{thm:k1ab}, in Section 4 we prove Theorem \ref{thm:w2k} and in Section 5 we prove Proposition \ref{thm:pkhat} and Theorem \ref{thm:p3p4p5hat}.  Finally, we present some concluding remarks in Section 6.
	
	{\bf Note.} After making this paper public, we learnt that Theorem~\ref{thm:k1ab} has also been proved by several other researchers  and Theorem~\ref{thm:w2k} has been proved by Methuku, Gr\'osz, Tompkins  (using a different proof). These works are unpublished.
	


\section{Preliminaries}
	In this section, we describe some preliminary tools that will be used in our proofs. The most important tool for proving Theorems \ref{thm:k1ab} and \ref{thm:w2k} is the triangle removal lemma~\cite{conlonFoxRemovalSurvey2013,foxNewRemoval2011,ruszaSzemerediTriple1976}. The specific form that we shall be using appears as Theorem 2.1 in \cite{conlonFoxRemovalSurvey2013}.
	\begin{lemma} [Ruzsa-Szemer\'edi~\cite{ruszaSzemerediTriple1976}]
		\label{lem:triangleRemoval}
		Suppose $\epsilon >0$. Let $\delta=\delta(\epsilon)$ be such that $\frac 1{\delta}$ is a tower of twos of height $684\log(\frac 1\epsilon)$. If $G$ is a graph on $n$ vertices with at least $\epsilon n^2$ edge-disjoint triangles, then $G$ contains at least $\delta n^3$ triangles.
	\end{lemma}
	
	 Recall that $t(G)$ is the number of triangles in $G$.
	\begin{lemma} [Nordhaus-Stewart~\cite{nordhausStewartTriangle1963}]
		\label{lem:nordhausStewart}
		For any graph $G$ on $n$ vertices,
		$$
		t(G)\ge \frac{e(G)}{3n}\cdot \left(4e(G)-n^2\right).$$
	\end{lemma}
	
	A $k$-uniform hypergraph $H=(V(H),E(H))$, consists of a vertex set $V(H)$ and edge set $E(H)\subseteq \binom {V(H)}k$. We write $e(H)=|E(H)|$. A subset $X\subseteq V(H)$ of vertices is an independent set if $e\not\subset X$ for any $e\in E(H)$. The independence number of $H$, denoted by $\alpha(H)$, is the largest size of an independent set in $H$.
	
	\begin{lemma} [Spencer~\cite{spencer-turan1972}]
		\label{lem:alphaH}
		Let $H=(V,E)$ be an $r$-uniform hypergraph on $n$ vertices, and let $d(H)=\frac{re(H)}{n}$ denote the average degree of $H$. Then, 
		$$
		\alpha(H)\ge\frac{r-1}{r}\cdot\left(\frac{n}{d(H)^{\frac1{r-1}}}\right).$$
	\end{lemma}
	
	Finally, we require a result that is a direct consequence of the proof of the Erd\H os-Gallai theorem \cite{erdosGallai1959} for cycles, which states that every graph with average degree at least $k$ contains a cycle of length at least $k+1$. Recall that a chord in a cycle is an edge between any two non-adjacent vertices of the cycle.
	
	\begin{lemma}[Erd\H os-Gallai \cite{erdosGallai1959}]
		\label{lem:erdosGallaiTheta}
		Let $k\ge 3$ be an integer. If $G$ is a graph of average degree at least $k$, then $G$ contains a cycle of length at least $k+1$ with a chord. In particular, $G$ also contains a path of length at least $k$.\hfill{$\Box$}
	\end{lemma}
	
	\section{Suspension of complete bipartite graphs}
	
	Our goal in this short section is to prove Theorem \ref{thm:k1ab} and give a lower bound on $\ex(n,K_3,K_{1,2,2})$ via Proposition \ref{prop:w4Lowerbd}. Given a graph $G$ and $X\subseteq V(G)$, let $N_G(X)=\bigcap_{v\in X}N_G(v)$ denote the common neighborhood of all vertices from $X$.
	
	\begin{proof}[Proof of Theorem \ref{thm:k1ab}]
		Recall that we are to show that for $1\le a\le b$, $\ex(n,K_3,K_{1,a,b})=o(n^{3-\frac1a})$.
		Let $n$ be sufficiently large, and $G$ be an $n$-vertex graph which is $K_{1,a,b}$-free. By the K\"{o}vari-S\'{o}s-Tur\'{a}n theorem \cite{kovariSosTuran1954} and (\ref{eq:introupper}), we know that there exists $c=c_b$ such that 
		\begin{equation} \label{eq:tupper}
		t(G)<c_bn^{3-\frac1a}. \end{equation}
		
		Now, suppose $\epsilon > 0$ is fixed. Assume that $n$ is sufficiently large, and that $G$ is a $K_{1,a,b}$-free graph  with $t(G) \ge \epsilon n^{3-\frac 1a}$. Construct an $a$-uniform hypergraph $H$ whose vertices are the triangles of $G$, and let $\{T_1,\ldots,T_a\}$  be an edge of $H$ if the triangles $T_1,\ldots,T_a$ all share a common edge in $G$. Then
		\[
		e(H) = \sum_{\{v_1,\ldots,v_a\}\in\binom{V(G)}a}e(N_G(\{v_1,\ldots,v_a\})).
		\]
		As $G$ is $K_{1,a,b}$-free, $N_G(\{v_1,\ldots,v_a\})$ has no vertex of degree at least $b$. This implies that $$e(N_G\{v_1, \ldots, v_a\})<\frac b2\cdot |N_G(\{v_1,\ldots,v_a\})|.$$ Consequently,
		\[
		e(H)<\sum_{\{v_1,\ldots,v_a\}\in\binom{[n]}a}\frac b2\cdot |N_G(v_1,\ldots,v_a)| = \frac b2 \sum_{v\in V(G)}\binom {\deg(v)} a < b\sum_{v\in V(G)} \deg(v)^a < b\cdot n^{a+1}.
		\]
		This implies that $d(H)<\frac {ab\cdot n^{a+1}}{t(G)}$. Using Lemma \ref{lem:alphaH},
		\[
		\alpha(H)> \frac{a-1}{a}\cdot \frac{t(G)}{\left(\frac {ab\cdot n^{a+1}}{t(G)}\right)^{\frac 1{a-1}}} = c\cdot {t(G)^{1+\frac1{a-1}}}\cdot {n^{-\frac{a+1}{a-1}}},
		\]
		where $c=(a-1)/(a(ab)^{1/(a-1)})$. Recalling that $t(G)\ge \epsilon n^{3-\frac 1a}$ and letting $\epsilon' = c\cdot \epsilon^{1+\frac 1{a-1}}$, we obtain
		\[
		\alpha(H)>\epsilon'n^{\frac a{a-1}\cdot \frac{3a-1}{a}}\cdot n^{-\frac{a+1}{a-1}} = \epsilon'n^2.
		\]
		
		Let $I$ be a maximum independent set of $H$. Create an auxiliary graph $H'$ with vertex set $I$, and join two vertices of $H'$ iff the triangles corresponding to them share an edge. Every triangle from $I$ can be adjacent to at most $3(a-1)$ other triangles from $I$. Therefore, $\deg_{H'}(i)<3a$ for every $i\in I$, and hence by Lemma \ref{lem:alphaH} $H'$ has an independent set of size at least $\frac{|I|}{6a}>\frac{\epsilon'n^2}{6a}$. The triangles corresponding to this independent set are edge-disjoint.
		Therefore $t(G) \ge \delta n^3$ where $\delta=\delta(\frac{\epsilon'}{6a})$ is obtained from Lemma \ref{lem:triangleRemoval}. However $t(G)<c_b n^{3-\frac 1a}$ by (\ref{eq:tupper}) and this implies that $\delta n^3 \le \frac{c_b}3\cdot n^{3-\frac 1a}$, a contradiction for sufficiently large $n$.
			\end{proof}
	
	Plugging in $a=b=2$ in Theorem \ref{thm:k1ab}, we get the bound $\ex(n,K_3,K_{1,2,2})=o(n^{5/2})$. We now describe a lower bound construction for $\ex(n,K_3,K_{1,2,2})$.
	
	\begin{prop}
		\label{prop:w4Lowerbd}
		When $n$ is a multiple of $4$,
		\[
		\ex(n,K_3,K_{1,2,2})\ge \frac{n^2}{4}.
		\]
	\end{prop}
	\begin{proof}[Proof of Proposition \ref{prop:w4Lowerbd}]
		Let $H_n=(A,B)$ be the complete bipartite graph with $|A|=|B|=\frac n2$, with additional edges in both the parts such that $H_n[A]$ and $H_n[B]$ are matchings of size $\frac n4$. Observe that the neighborhood of every vertex in $H_n$ consists of $\frac n4$ edge-disjoint triangles sharing a common vertex, and hence is $C_4$-free. Thus $H_n$ is $K_{1,2,2}$-free. On the other hand, every triangle of $H_n$ either has an edge inside $A$ or an edge inside $B$, implying
		\[
		t(H_n) = e(A)\cdot |B| + e(B)\cdot |A| = 2\cdot \frac n4\cdot \frac n2 = \frac{n^2}{4}.
		\]
		This proves that whenever $4\mid n$, $\ex(n,K_3,K_{1,2,2})\ge n^2/4$.
	\end{proof}

	\section{Suspension of even cycles}
	Our goal in this section is to prove Theorem \ref{thm:w2k}.
	Before proceeding with the proof, we prove Lemma \ref{lem:pathsLemma}  which gives an upper bound on the number of paths of length $k$ in a $C_{2k}$-free graph. The main idea behind the lemma is the technique used in \cite{pikhurkoEvenCycle2012,verstraeteUnavoidable2005,verstraeteCycleSurvey2016} to prove upper bounds on $\ex(n,C_{2k})$ by analyzing the breadth-first search tree from any vertex.
	
	Given a graph $G$ and a vertex $r\in V(G)$, a breadth-first search tree $T$ of $G$ rooted at $r$ is constructed as follows. Let $L_0=\{r\}$. For $i\ge 1$, let $L_i\subseteq V(G)$ be the set of all vertices in $V(G)$ which are at distance $i$ from vertex $r$. The vertex subset $L_i$ is called the $i$'th level of $T$. The tree $T$ consists of vertex set $V(G)$ and edge set a subset of the edges of $G$ between levels $L_i$ and $L_{i+1}$, $i\ge 0$, such that every vertex in $L_{i+1}$ is incident to exactly one vertex in $L_i$.
	
	For $i\ge 0$, let $G[L_i]$ be the subgraph of $G$ induced by $L_i$, and let $G[L_i,L_{i+1}]$ be the bipartite subgraph of $G$ with parts $(L_i,L_{i+1})$ and edges exactly those of $G$ with one endpoint in $L_i$ and another in $L_{i+1}$.
	
	We now quote Lemma 3.5 from \cite{verstraeteCycleSurvey2016} in the form that we shall be using.
	
	\begin{lemma} [Verstra\"ete~\cite{verstraeteCycleSurvey2016}]
		\label{lem:bfsDeg}
		Let $T$ be a breadth-first search tree in a graph $G$, with levels $L_0, L_1, \ldots$. Let $\Theta_\ell$ denote the set of cycles of length $\ell$ with a chord. Then,
		
		(a) If $G[L_i]$ contains an element of $\Theta_\ell$, then there exists $m\le i$ such that $G$ contains cycles $C_{2m+1}$ ,$C_{2m+2}$, $\ldots$, $C_{2m+\ell-1}$.
		
		(b) If $G[L_i,L_{i+1}]$ contains an element of $\Theta_\ell$, then there exists $m\le i$ such that $G$ contains cycles $C_{2m+2}$, $C_{2m+4}$, $\ldots$, $C_{2m+\ell'}$, where $\ell'$ is the largest even integer less than $\ell$.
	\end{lemma}
	Let $p_k(G)$ denote the number of paths of length $k$ in a graph $G$. Here each subgraph isomorphic to $P_k$ is counted twice, once for each ordering of its vertices along the path.
	
	\begin{lemma}
		\label{lem:pathsLemma}
		Let $k \ge 2$ be an integer, $0<\epsilon<1$ and $n>(20k/\epsilon)^k$. Let $F$ be a $C_{2k}$-free graph on $n$ vertices with minimum degree at least $\epsilon n^{1/k}$. Then
		\[
		p_k(F)\le \left(\frac{2k}{\epsilon}\right)^{(k-1)k}n^2.
		\]
	\end{lemma}
	
	\begin{proof}[Proof of Lemma \ref{lem:pathsLemma}]
		We first prove that  $F$ has  bounded maximum degree. Suppose for contradiction that there exists $v \in V(F)$ with 
		\[
		\deg(v)\ge \left(\frac{2k}{\epsilon}\right)^{k-1}\cdot n^{1/k}.
		\]
		
		Consider the breadth-first search tree of $F$ starting at $v$. For $i\ge 0$, let $L_i$ be the $i$th level of this breadth-first search tree. By assumption, $|L_1|\ge \left(\frac{2k}{\epsilon}\right)^{k-1}\cdot n^{1/k}$. Denote by $e(L_i,L_{i+1})$ the number of edges in $G[L_i,L_{i+1}]$. Let us prove  that for every $1\le i<k$,
		\begin{equation} \label{eq:lili}
		e(L_i)\le (k-1)|L_i| \quad \mbox{ and } \quad e(L_i,L_{i+1})\le (k-1)(|L_i|+|L_{i+1}|).
		\end{equation}
		Indeed, if $e(L_i)>(k-1)|L_i|$, then by Lemma \ref{lem:erdosGallaiTheta}, $F[L_i]$ contains a cycle of length $\ell$ with a chord, where $\ell\ge 2k-1$. Now, we apply Lemma \ref{lem:bfsDeg}{ (a)} to obtain an integer $m\le i$ such that $F$ contains cycles of lengths $2m+1, 2m+2,\ldots, 2m+\ell -1$. Since $\ell \ge 2k-1$ and $1\le m<k$, we have $2m+1\le 2k \le 2m+\ell -1$. Then $F$ contains a $C_{2k}$, contradiction.
		Similarly, if $e(L_i,L_{i+1})>(k-1)(|L_{i}|+|L_{i+1}|)$, then Lemma \ref{lem:erdosGallaiTheta} gives us a cycle of length $\ell$ with a chord in $F[L_i,L_{i+1}]$ where $\ell\ge 2k-1$. Then Lemma \ref{lem:bfsDeg}{ (b)} gives an integer $m\le i$ such that $F$ contains cycles of lengths $2m+2,2m+4,\ldots, 2m+\ell$. As $\ell\ge 2k-1$ and $1\le m<k$, we have $2m+2\le 2k \le 2m+2k-1$. This implies that $F$ contains a $C_{2k}$, again a contradiction.
		
		\begin{clm}
			\label{clm:Li+1Limult}
			For every $i\ge 0$, 
			\begin{equation}
			\label{eq:Li+1Limult}
			|L_{i+1}|\ge \frac{\epsilon n^{1/k}}{2k}\cdot |L_i|.
			\end{equation}
		\end{clm}
		\noindent {\it Proof of Claim \ref{clm:Li+1Limult}.}
		We use induction on $i$. Note that $|L_0|=1$ and 
		$$|L_1|=\deg(v)\ge \left(\frac{2k}{\epsilon}\right)^{k-1}\cdot n^{1/k} >\frac{\epsilon n^{1/k}}{2k} =\frac{\epsilon n^{1/k}}{2k}|L_0|.$$ Moreover, for $i\ge 1$ and any vertex $u\in L_i$, $N_G(u)\subseteq L_{i-1}\cup L_i\cup L_{i+1}$. Thus, (\ref{eq:lili}) implies 
		\begin{equation}
		\begin{split}
		k(|L_i|+|L_{i+1}|) + 2k|L_i| + k(|L_i|+|L_{i-1}|) &> e(L_i,L_{i+1}) + 2e(L_i) + e(L_i,L_{i-1}) \notag \\
		& =\sum_{u\in L_i}\deg(u) \notag \\ &\ge \epsilon n^{1/k}\cdot |L_i|. \notag \\
		\end{split}
		\end{equation}
		Consequently, 
		\begin{equation}
		\label{eq:Li+1Li}
		|L_{i+1}|>\left(\frac{\epsilon n^{1/k}}{k}-4\right)\cdot|L_i|-|L_{i-1}|.
		\end{equation}
		By the induction hypothesis we may assume that $|L_{i-1}|\le \frac{2k}{\epsilon n^{1/k}}\cdot |L_{i}|$. Thus, (\ref{eq:Li+1Li}) implies
		\[
		|L_{i+1}|>\left(\frac{\epsilon n^{1/k}}{k}-4-\frac{2k}{\epsilon n^{1/k}}\right)\cdot|L_i| > \frac{\epsilon n^{1/k}}{2k}\cdot |L_i|
		\]
		since $n>(20k/\epsilon)^k$. This finishes the proof of Claim \ref{clm:Li+1Limult}.
		
		Now, by applying Claim \ref{clm:Li+1Limult} iteratively, we obtain
		\[
		|L_k|\ge \left(\frac{\epsilon n^{1/k}}{2k}\right)^{k-1}\cdot |L_1|\ge \left(\frac{\epsilon n^{1/k}}{2k}\right)^{k-1}\cdot \left(\frac{2k}{\epsilon}\right)^{k-1}\cdot n^{1/k}=n,
		\]
		a contradiction. Thus, $\deg(u)\le  \left(\frac{2k}{\epsilon}\right)^{k-1}\cdot n^{1/k}$ for every $u\in V(F)$. Therefore, if $\Delta(F)$ is the maximum degree of $F$,
		\[
		p_k(F)\le n\cdot \Delta(F)^k \le \left(\frac{2k}{\epsilon}\right)^{k(k-1)}\cdot n^2,
		\]
		as desired.
	\end{proof}
	
 For a graph $G$ and edge $uv\in E(G)$, the codegree of $uv$ is  $\deg_G(u,v) = |N_G(\{u,v\})|$. 
	
	\begin{proof}[Proof of Theorem \ref{thm:w2k}]
		Fix $\epsilon >0$ and let $n$ be sufficiently large.
		Suppose $G=(V,E)$ is a graph on $n$ vertices satisfying $t(G)\ge \epsilon n^{2+\frac 1k}$. We wish to show that $G$ contains a copy of $\widehat{C}_{2k}$. Suppose, on the contrary, that $G$ is $\widehat{C}_{2k}$-free. 
			First, we iteratively delete edges of $G$ with codegree less than $\frac{\epsilon n^{1/k}}{10}$ in the current graph, until there are no such edges left. Since we delete fewer than $e(G)\cdot \frac{\epsilon n^{1/k}}{10}$ triangles, we are left with a graph $G'$ satisfying
		\[
		t(G') >t(G)-e(G)\cdot \frac{\epsilon n^{1/k}}{10} > \epsilon n^{2+\frac 1k}-\frac{\epsilon n^{2+\frac 1k}}{10} \ge \frac{9\epsilon}{10}\cdot n^{2+\frac 1k},
		\]
		and $\deg_{G'}(u,v)\ge \frac{\epsilon n^{1/k}}{10}$ for every  $uv\in E(G')$.
		
		Next, create an auxiliary graph $H$ whose vertices are the triangles of $G'$, and two vertices of $H$ are adjacent iff their corresponding triangles share a common edge. By Lemma \ref{lem:alphaH}, $\alpha(H)>\frac{t(G')}{2d(H)}$. Let \[\gamma := \frac{(\epsilon/20)^{k^2}}{k^{(k-1)k}}>0.\]
		If $\frac{t(G')}{2d(H)}>\frac\gamma2 n^2$, then this gives us $\frac\gamma2 n^2$ edge-disjoint triangles in $G$, and  this implies that $t(G) > \delta n^3$  where $\delta=\delta(\frac \gamma2)$ from Lemma~\ref{lem:triangleRemoval}. However, we  also have $t(G) < c_k n^{2+\frac 1k}$ for some constant $c_k$ by 	(\ref{eq:introupper}) and this is a contradiction since $n$ is sufficiently large. Therefore, we may assume
		$\frac{t(G')}{d(H)}\le \gamma n^2$ and this implies 
		\[ d(H)\ge \frac{t(G')}{\gamma n^2} \ge \frac{9\epsilon}{10\gamma}\cdot n^{1/k}
		\]
		 and hence
		\[e(H)\ge \frac{9\epsilon n^{1/k}}{20\gamma}\cdot t(G')\ge \frac{81\epsilon}{200\gamma}\cdot n^{2+\frac 2k}.\]
		
		Let us now bound $X$, the number of copies of $\widehat{P}_k$ in $G'$ in two different ways. For every $v\in V(G')$, let $G'_v$ denote the subgraph of $G'$ induced by $N_{G'}(v)$. Let $\delta(F)$ denote the minimum degree of $F$ for any graph $F$. By the assumption on the minimum codegree of edges in $G'$, $\delta(G'_v)\ge \frac{\epsilon n^{1/k}}{10}$. Hence applying Lemma \ref{lem:pathsLemma} with $\epsilon$ replaced by $\frac \epsilon{10}$,
		\begin{equation}
		\label{eq:Xleft}
		\begin{aligned}
		X\le \sum_{v\in V(G')}p_k(G'_v) &\le n\cdot \left(\frac{20k}{\epsilon}\right)^{k(k-1)}\cdot n^2\\& = \left(\frac{20k}{\epsilon}\right)^{k(k-1)} \cdot n^3.
		\end{aligned}
		\end{equation}
		On the other hand, we can first fix two adjacent triangles in $G'$ and then keep growing it to a $\widehat{P}_k$ by using the minimum codegree condition of $G'$. Since $\delta(H)\ge \frac{\epsilon n^{1/k}}{10}$, this implies that for large enough $n$,
		\begin{equation}
		\label{eq:Xright}
		\begin{aligned}
		X&\ge \frac 12 e(H)\cdot (\delta(H)-2)\cdot (\delta(H)-3)\cdots (\delta(H)-k+1)\\
		 &\ge \frac 12 e(H)\cdot (\delta(H)-k)^{k-2}\\
		 &\ge \frac{81\epsilon}{400\gamma}\cdot n^{2+\frac 2k}\cdot \left(\frac{\epsilon n^{1/k}}{20}\right)^{k-2}\\
		 &= \frac{81\epsilon^{k-1}}{20^k\cdot \gamma}\cdot  n^3.
		\end{aligned}
		\end{equation}
		The factor $\frac 12$ in  (\ref{eq:Xright}) is to balance out over-counting the same $\widehat{P}_k$ from its two ends. Comparing (\ref{eq:Xleft}) and (\ref{eq:Xright}), we obtain
		\[
		\left(\frac{20k}{\epsilon}\right)^{k(k-1)}\ge \frac{81\epsilon^{k-1}}{20^k\cdot \gamma},
		\]
		implying
		\[
		\gamma \ge \frac{81\epsilon^{k^2}}{20^{k^2}\cdot k^{k(k-1)}}=81\cdot \frac{(\epsilon/20)^{k^2}}{k^{k(k-1)}}=81\gamma,
		\]
		a contradiction. This completes the proof of Theorem \ref{thm:w2k}.
	\end{proof}


\section{Suspension of paths}
	In this section, we prove Proposition~\ref{thm:pkhat} and Theorem \ref{thm:p3p4p5hat}. 
	
	\subsection{Proof of Proposition~\ref{thm:pkhat}}
		First, we show the upper bound in (\ref{eq:thm-pkhat-upper}). Let $k\ge 3$ be fixed, and let $G$ be a graph on $n$ vertices which is $\widehat{P}_k$-free. We need to show that $t(G)\le \frac{(k-1)}{12}\cdot n^2 + \frac{(k-1)^2}{12}\cdot n$. 
		\smallskip
		
		Note that the neighborhood of every vertex $v\in V(G)$ is $P_k$-free. Thus by Lemma \ref{lem:erdosGallaiTheta}, the average degree of the subgraph of $G$ induced by $N_G(v)$ is at most $k-1$. Hence,
		\[
		e(N_G(v))\le \frac{k-1}{2}\cdot \deg_G(v).
		\]
		Summing up this inequality over all vertices $v\in V(G)$,
		\[
		3t(G)=\sum_{v\in V(G)}e(N_G(v))\le \frac{k-1}{2}\cdot 2e(G) = (k-1)e(G),
		\]
		giving us
		\begin{equation}
		\label{eq:triangleEdge}
		t(G)\le \frac{k-1}3\cdot e(G).
		\end{equation}
		
		This, in conjunction with Lemma \ref{lem:nordhausStewart}, gives us
		\[
		\frac{k-1}{3}\cdot e(G)\ge t(G)\ge \frac{e(G)}{3n}\cdot(4e(G)-n^2),
		\]
		which simplifies to
		\[e(G)\le \frac{n^2}{4}+\frac{(k-1)n}{4}.\]
		The conclusion of the upper bound follows from plugging this inequality back into (\ref{eq:triangleEdge}).
		
		\medskip
		Now we prove the lower bound in 
		(\ref{eq:thm-pkhat-upper}). Let $n$ be a multiple of $4\left\lfloor \frac{k-1}{2}\right\rfloor$. We shall construct a $\widehat{P}_k$-free graph $F_{n,k}$ on $n$ vertices with $t(F_{n,k})\ge \left\lfloor\frac{k-1}{2}\right\rfloor\cdot\frac{n^2}{8}$.
		
		Let $F_{n,k}=(A,B)$ be the complete bipartite graph with parts $A, B$ with $|A|=|B|=\frac n2$, with additional edges in  $A$ such that $F_{n,k}[A]$ is a disjoint union of $K_{\left\lfloor\frac{k-1}{2}\right\rfloor,\left\lfloor\frac{k-1}{2}\right\rfloor}$. Then,
		\[
		\textstyle e(A) = \left\lfloor\frac{k-1}{2}\right\rfloor^2\cdot \frac{n}{4\left\lfloor\frac{k-1}2\right\rfloor}=\left\lfloor\frac{k-1}2\right\rfloor \cdot\frac n4.
		\]
		Every triangle of $F_{n,k}$ consists of an edge from $F_{n,k}[A]$ and a vertex from $B$. Hence,
		\[
		\textstyle t(F_{n,k}) = \left\lfloor\frac{k-1}2\right\rfloor \cdot\frac n4\cdot \frac n2 = \left\lfloor\frac{k-1}2\right\rfloor \cdot\frac{n^2}{8}.
		\]
		Further, $F_{n,k}$ is $\widehat{P}_k$-free, since the neighborhood of every vertex in $B$ is a disjoint union of $K_{\left\lfloor\frac{k-1}{2}\right\rfloor,\left\lfloor\frac{k-1}{2}\right\rfloor}$, and the neighborhood of every vertex in $A$ is isomorphic to $K_{\left\lfloor\frac{k-1}{2}\right\rfloor,\frac n2}$.	\hfill{$\Box$}
	
\subsection{Proof of Theorem \ref{thm:p3p4p5hat}}	
  We will use some ideas from \cite{methukuK3C52019}, and define the concepts of \emph{triangle-connectivity} and \emph{blocks}. In what follows, a triangle $T$ in a graph $G$ is a set of three edges $\{ab,bc,ca\}$ that form a $K_3$ in $G$. Subsequently, we shall denote such a triangle simply as $abc$.
	
	\begin{defn}[Triangle-connectivity]
		 Given a graph $G$ and two distinct edges $e, e'\in E(G)$,  say that $e$ and $e'$ are \emph{triangle-connected} if there is a sequence of triangles $\{T_1,\ldots,T_k\}$ of $G$, such that $e\in T_1$, $e'\in T_k$, and $T_i$ and $T_{i+1}$ share a common edge for every $1\le i\le k-1$. A subgraph $H\subseteq G$ is \emph{triangle-connected} if $e$ and $e'$ are triangle-connected for every two distinct $e, e'\in E(H)$.
	\end{defn}
	
	It is straightforward to check that triangle-connectivity is an equivalence relation on $E(G)$ (assuming reflexivity as part of the definition).
	
	\begin{defn}[Triangle block]
		A \emph{triangle block}, or simply a \emph{block} in a graph $G$ is a subgraph $H$ whose edges form an equivalence class of the \emph{triangle-connectivity} relation on $E(G)$.
	\end{defn}
	
	In other words, a subgraph $H\subseteq G$ is a triangle block if it is edge-maximally triangle-connected. By definition, the triangle blocks of a graph $G$ are edge-disjoint. 
	\subsubsection{Proof of Theorem \ref{thm:p3p4p5hat} for $\mathbf{ k=3}$.}
		Suppose $G$ is a graph on $n$ vertices which is $\widehat{P}_3$-free. We will prove using induction on $n$, that 
		\begin{equation}
		\label{eq:p3hat-induction}
		t(G)< \frac{n^2}{8}+3n.
		\end{equation}
		
		This inequality is true for $n=3$ as $t(G)\le 1$ for any graph $G$ on $3$ vertices. Now, fix an $n>3$ and a graph $G$ on $n$ vertices which is $\widehat{P}_3$-free. Assume that (\ref{eq:p3hat-induction}) holds for $\widehat{P}_3$-free graphs on less than $n$ vertices.
		
		We may assume without loss of generality that every edge of $G$ lies in a triangle, otherwise we may delete it from $G$ without changing $t(G)$. For a vertex $v\in V(G)$, let $t(v) = t_G(v)$ denote the number of triangles in $G$ containing $v$. By definition, $t_G(v) = e(N_G(v))$.
		
		We first prove that if $G$ has a copy of $K_4$, then $t(G)<\frac{n^2}8+3n$, hence completing the induction step.
		
		Suppose $G$ has a copy of $K_4$ with vertices labeled $a_1,a_2,a_3,a_4$. Let $X = V(G)\setminus \{a_1,a_2,a_3,a_4\}$, and $A_i = N_G(a_i)\cap X$, $i=1,\ldots,4$. If $x\in A_1\cap A_2$, then we can find a $\widehat{P}_3$ formed by $a_1,x,a_2,a_3,a_4$ in the neighborhood of $a_1$. Thus, $A_1\cap A_2=\varnothing$, and by symmetry the $A_i$'s are mutually disjoint. Hence, $|A_1|+|A_2|+|A_3|+|A_4|\le |X| = n-4$. This implies that one of the $A_i$'s has size $\le \frac{n-4}{4}$. Using Lemma \ref{lem:erdosGallaiTheta} in the neighborhood of $a_i$,
		\[
		t(a_i) = 3 + e(A_i) \le 3 + |A_i| \le 3 + \frac{n-4}{4} = \frac{n+8}{4}.
		\]
		Now let $G'=G-a_i$. As $G$ was $\widehat{P}_3$-free, so is $G'$. Hence by the induction hypothesis,
		\[
		t(G')< \frac{(n-1)^2}{8}+3(n-1).
		\]
		This implies,
		\[
		t(G)=t(G')+t(a_i) < \frac{(n-1)^2}8+3(n-1)+\frac{n+8}4 < \frac{n^2}8+3n,
		\]
		as desired.		We may now assume that $G$ is $K_4$-free. 

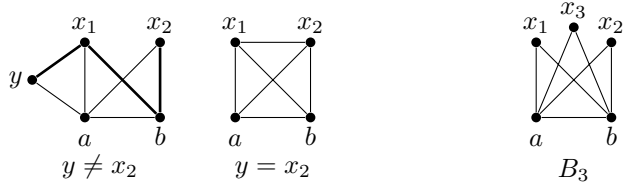
\begin{figure}[h]
			\centering
			\begin{tikzpicture}
			\tikzstyle{vertex}=[circle,fill=black,minimum size=2pt,inner sep=1.3pt]
			\node[vertex] (a) at (0,0){};
			\node[vertex] (b) at (1,0){};
			\node[vertex] (x1) at (0,1){};
			\node[vertex] (x2) at (1,1){};
			\node[vertex] (y) at (-0.7,0.5){};
			\draw (a)--(x1)--(b)--(x2)--(a)--(b);
			\draw [line width = 1pt] (y)--(x1)--(b)--(x2);
			\draw (x1)--(y)--(a);
			\draw (0,-0.1) node [below] {$a$};
			\draw (b) node [below] {$b$};
			\draw (x1) node [above] {$x_1$};
			\draw (x2) node [above] {$x_2$};
			\draw (y) node [left] {$y$};
			\draw (0.2,-0.65)  node {$y\neq x_2$};
			
			\node[vertex] (ap) at (2,0){};
			\node[vertex] (bp) at (3,0){};
			\node[vertex] (x1p) at (2,1){};
			\node[vertex] (x2p) at (3,1){};
			\draw (ap)--(x1p)--(bp)--(x2p)--(ap)--(bp);
			\draw (x1p)--(x2p);
			\draw (2,-0.1) node [below] {$a$};
			\draw (bp) node [below] {$b$};
			\draw (x1p) node [above] {$x_1$};
			\draw (x2p) node [above] {$x_2$};
			\draw (2.5,-0.7)  node {$y=x_2$};
			
			\node[vertex] (app) at (6,0){};
			\node[vertex] (bpp) at (7,0){};
			\node[vertex] (xpp) at (6,1){};
			\node[vertex] (ypp) at (7,1){};
			\node[vertex] (zpp) at (6.5,1.2){};
			\draw (app)--(bpp)--(xpp)--(app)--(ypp)--(bpp)--(zpp)--(app);
			\draw (6,-0.1) node [below] {$a$};
			\draw (bpp) node [below] {$b$};
			\draw (xpp) node [above] {$x_1$};
			\draw (ypp) node [above] {$x_2$};
			\draw (zpp) node [above] {$x_3$};
			\draw (6.5,-0.7)  node {$B_3$};
			\end{tikzpicture}
			\caption{(left): third triangle on $ax_1$, (right): third triangle on $ab$\label{fig:blockp3}}
		\end{figure}

		 Let $B_s$ denote the \emph{book graph} on $s+2$ vertices, consisting of $s$ triangles all sharing a common edge.
		\begin{clm}
			\label{clm:blockp3}
			Every triangle block of $G$ is isomorphic to $B_s$ for some $s\ge 1$.
		\end{clm}
		\noindent {\it Proof of Claim \ref{clm:blockp3}.}		
		Let $H\subseteq G$ be an arbitrary triangle block. If $H$ contains only one or two triangles, it is isomorphic to $B_1$ or $B_2$. Suppose $H$ contains at least three triangles. Let two of them be $abx_1$ and $abx_2$ (Figure \ref{fig:blockp3}). If another triangle is of the form $ax_1y$ for some $y\in V(H)$, then there are two possible cases. If $y\neq x_2$, then $N_H(a)$ contains the 3-path $x_2bx_1y$. Otherwise, if $y=x_2$, then the vertices $a,b,x_1,x_2$ create a $K_4$. Similarly, no triangle contains any of the edges $bx_1, ax_2, bx_2$. Therefore all triangles in $H$ contain  $ab$ and $H\cong B_s$ for some $s\ge 1$. 
		
		\medskip
		Claim \ref{clm:blockp3} implies that $G$ comprises  $r$ edge-disjoint blocks isomorphic to books for some $r\ge 1$. Let the blocks of $G$ be isomorphic to $B_{s_1},\ldots,B_{s_r}$, where $s_1,\ldots, s_r\ge 1$. Then,
		\[
		t(G)=s_1+\cdots+s_r \quad \mbox{ and } \quad e(G)=2(s_1+\cdots+s_r)+r=2t(G)+r.
		\]
		Hence, $t(G)<e(G)/2$.		Finally, we apply Lemma \ref{lem:nordhausStewart} on $G$ to obtain
		\[
		\frac{e(G)}2>t(G)\ge \frac{e(G)}{3n}\cdot(4e(G)-n^2),
		\]
		implying
		\[
		e(G)<\frac{n^2}{4}+\frac{3n}{8}.
		\]
		Therefore $t(G)<\frac{n^2}8+\frac{3n}{16}<\frac{n^2}8+3n$, completing the induction step.
		\hfill{$\Box$}
		
	\subsubsection{Proof of Theorem \ref{thm:p3p4p5hat} for $\mathbf{k=4}$.}
		Suppose $\epsilon>0$ and $n$ is sufficiently large. Let $G$ be any graph on $n$ vertices which is $\widehat{P}_4$-free, such that \[t(G)\ge\frac{n^2}{8}+13\epsilon n^2.\]
		
		
		The very first step of the proof is to remove copies of $K_4$ and $K_{1,2,2}$ from $G$ while still maintaining $t(G)\ge \frac{n^2}8+\epsilon n^2$. We achieve this by means of the triangle removal lemma. First, we make an observation which follows immediately from Lemma~\ref{lem:triangleRemoval}, (\ref{eq:introupper}) and Lemma~\ref{lem:erdosGallaiTheta} for large $n$.
		\begin{equation}
		\label{eq:notManyEdgeDisjointTriangles}
		\mbox{$G$ cannot have }\epsilon n^2\mbox{ edge-disjoint triangles.}
		\end{equation}

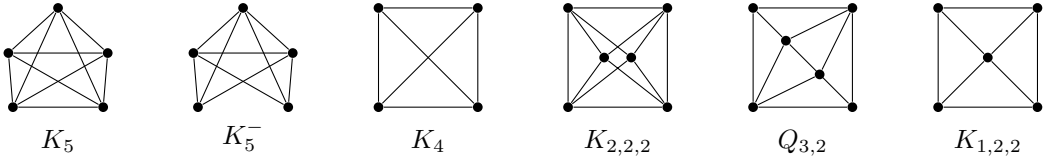
\begin{figure}[h]
			\centering
			\begin{tikzpicture}[scale = 0.6]
			\tikzstyle{vertex}=[circle,fill=black,minimum size=2pt,inner sep=1.3pt]
			\begin{scope}[shift={(-4.1,0)}]
			\node[vertex] (u) at (4.1,1) {};
			\node[vertex] (a1) at (3,0) {};
			\node[vertex] (b1) at (5.2,0) {};
			\node[vertex] (c1) at (5.1,-1.2) {};
			\node[vertex] (d1) at (3.1,-1.2) {};
			\draw (a1)--(b1)--(c1);
			\draw (d1)--(a1)--(u)--(b1);
			\draw (a1)--(c1);
			\draw (b1)--(d1);
			\draw (c1)--(u)--(d1);
			\draw (4.1,-1.3) node [below] {$K^-_5$};
			\end{scope}
			
			\begin{scope}[shift={(-8.2,0)}]
			\node[vertex] (u) at (4.1,1) {};
			\node[vertex] (a1) at (3,0) {};
			\node[vertex] (b1) at (5.2,0) {};
			\node[vertex] (c1) at (5.1,-1.2) {};
			\node[vertex] (d1) at (3.1,-1.2) {};
			\draw (a1)--(b1)--(c1)--(d1)--(a1)--(u)--(b1);
			\draw (a1)--(c1);
			\draw (b1)--(d1);
			\draw (c1)--(u)--(d1);
			\draw (4.1,-1.5) node [below] {$K_5$};
			\end{scope}
			
			\node[vertex] (a1) at (3,1) {};
			\node[vertex] (b1) at (5.2,1) {};
			\node[vertex] (c1) at (5.2,-1.2) {};
			\node[vertex] (d1) at (3,-1.2) {};
			\draw (a1)--(b1)--(c1)--(d1)--(a1);
			\draw (a1)--(c1);
			\draw (b1)--(d1);
			\draw (4.1,-1.5) node [below] {$K_4$};
			
			\node[vertex] (a3) at (7.2,1) {};
			\node[vertex] (b3) at (9.4,1) {};
			\node[vertex] (c3) at (9.4,-1.2) {};
			\node[vertex] (d3) at (7.2,-1.2) {};
			\node[vertex] (u3) at (8.0,-0.1) {};
			\node[vertex] (v3) at (8.6,-0.1) {};
			\draw (a3)--(b3)--(c3)--(d3)--(a3);
			\draw (a3)--(u3)--(b3)--(v3)--(c3);
			\draw (c3)--(u3)--(d3)--(v3)--(a3);
			\draw (8.3,-1.5) node[below] {$K_{2,2,2}$};
			
			\begin{scope}[shift={(4.1,0)}]
			\node[vertex] (a3) at (7.2,1) {};
			\node[vertex] (b3) at (9.4,1) {};
			\node[vertex] (c3) at (9.4,-1.2) {};
			\node[vertex] (d3) at (7.2,-1.2) {};
			\node[vertex] (u3) at ($(a3)!0.33!(c3)$) {};
			\node[vertex] (v3) at ($(a3)!0.67!(c3)$) {};
			\draw (a3)--(b3)--(c3)--(d3)--(a3);
			\draw (a3)--(u3)--(b3);
			\draw (c3)--(u3)--(d3);
			\draw (d3)--(v3)--(b3);
			\draw (8.3,-1.5) node[below] {$Q_{3,2}$};
			\end{scope}
			
			\begin{scope}[shift={(8.2,0)}]
			\node[vertex] (a3) at (7.2,1) {};
			\node[vertex] (b3) at (9.4,1) {};
			\node[vertex] (c3) at (9.4,-1.2) {};
			\node[vertex] (d3) at (7.2,-1.2) {};
			\node[vertex] (u3) at (8.3,-0.1) {};
			\draw (a3)--(b3)--(c3)--(d3)--(a3);
			\draw (a3)--(u3)--(b3);
			\draw (c3)--(u3)--(d3);
			\draw (8.3,-1.5) node[below] {$K_{1,2,2}$};
			\end{scope}
			\end{tikzpicture}
			\caption{\label{fig:k4w42w4}The graphs $K_5$, $K_5^-$, $K_4$, $K_{2,2,2}$, $Q_{3,2}$ and $K_{1,2,2}$.}
		\end{figure}

		Without loss of generality, we may assume that every edge of $G$ is contained in a triangle. We shall use (\ref{eq:notManyEdgeDisjointTriangles}) to remove all copies of the following six graphs in this order: $K_5$; $K_5^-$; $K_4$; $K_{2,2,2}$; $Q_{3,2}=\overbar{K_2}\vee P_3$; and $K_{1,2,2}$ (Figure \ref{fig:k4w42w4}).
		
	\begin{itemize}
		\item {\bf Step 1:} Cleaning $K_5$'s.
		
		If $G$ contains a $K_5$ with vertices $a_1,a_2,a_3,a_4,a_5$, then it has to be a block by itself. This is because if there is a vertex $x\neq a_i$ with $xa_1,xa_2\in E(G)$, then $N_G(a_1)$  contains the path $a_5a_4a_3a_2x$,  contradiction. Hence, all the $K_5$'s in $G$ are edge-disjoint.
		
		If $G$ has more than $\epsilon n^2$ copies of $K_5$, then by taking one triangle from each $K_5$, we get $\epsilon n^2$  edge-disjoint triangles in $G$, contradicting (\ref{eq:notManyEdgeDisjointTriangles}). Therefore, $G$ has at most $\epsilon n^2$ copies of $K_5$.
		
		We now delete one edge from each copy of $K_5$ in $G$, and lose at most $3\epsilon n^2$ triangles from $G$. So, we may assume $t(G)\ge \frac{n^2}{8}+10\epsilon n^2$, and $G$ is $\{\widehat{P}_4,K_5\}$-free.
		
		\item {\bf Step 2:} Cleaning $K_5^-$'s.
		
		Suppose $G$ contains a $K_5^-$. Observe that if we have a new vertex $x\neq a_i$ which is adjacent to two endpoints of any edge of this $K_5^-$, it would create a copy of $\widehat{P}_4$ (see Figure \ref{fig:k5-k4edgedisjoint} (left)). Thus, the only way two $K_5^-$'s can intersect in an edge is if they share the same five vertices. This would give us a $K_5$ in $G$, a contradiction. Therefore, the copies of $K_5^-$ are all edge-disjoint.
		
		Hence, if $G$ has more than $\epsilon n^2$ copies of $K_5^-$, we  again obtain at least $\epsilon n^2$ edge-disjoint triangles in $G$, contradicting (\ref{eq:notManyEdgeDisjointTriangles}). So $G$ has at most $\epsilon n^2$ copies of $K_5^-$.
		
		Deleting one edge from each copy of $K_5^-$ in $G$, we lose at most $3\epsilon n^2$ triangles in the process. After deletion, we still have $t(G)\ge \frac{n^2}{8}+7\epsilon n^2$, and we can further assume that $G$ is $\{\widehat{P}_4,K_5^-\}$-free.
		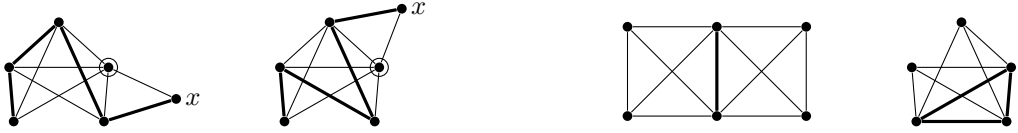
\begin{figure}[h]
			\centering
			\begin{tikzpicture}[scale = 0.6]
			\tikzstyle{vertex}=[circle,fill=black,minimum size=2pt,inner sep=1.3pt]
			\begin{scope}
			\node[vertex] (u) at (4.1,1) {};
			\node[vertex] (a1) at (3,0) {};
			\node[vertex] (b1) at (5.2,0) {};
			\node[vertex] (c1) at (5.1,-1.2) {};
			\node[vertex] (d1) at (3.1,-1.2) {};
			\node[vertex] (v) at (6.7,-0.7) {};
			\draw (a1)--(b1)--(c1);
			\draw (d1)--(a1)--(u)--(b1);
			\draw (a1)--(c1);
			\draw (b1)--(d1);
			\draw (c1)--(u)--(d1);
			\draw (c1)--(v)--(b1);
			\draw(b1)circle(0.2);
			\draw [line width = 1.2pt] (v)--(c1)--(u)--(a1)--(d1);
			\draw (v) node[right] {$x$};
			\end{scope}
			\begin{scope}[shift={(6,0)}]
			\node[vertex] (u) at (4.1,1) {};
			\node[vertex] (a1) at (3,0) {};
			\node[vertex] (b1) at (5.2,0) {};
			\node[vertex] (c1) at (5.1,-1.2) {};
			\node[vertex] (d1) at (3.1,-1.2) {};
			\node[vertex] (v) at (5.7,1.3) {};
			\draw (a1)--(b1)--(c1);
			\draw (d1)--(a1)--(u)--(b1);
			\draw (a1)--(c1);
			\draw (b1)--(d1);
			\draw (c1)--(u)--(d1);
			\draw (b1)--(v)--(u);
			\draw(b1)circle(0.2);
			\draw [line width = 1.2pt] (v)--(u)--(c1)--(a1)--(d1);
			\draw (v) node[above,right] {$x$};
			\end{scope}
			\begin{scope}[shift={(20,0)}]
			\node[vertex] (u) at (4.1,1) {};
			\node[vertex] (a1) at (3,0) {};
			\node[vertex] (b1) at (5.2,0) {};
			\node[vertex] (c1) at (5.1,-1.2) {};
			\node[vertex] (d1) at (3.1,-1.2) {};
			\draw (a1)--(b1)--(c1)--(d1)--(a1);
			\draw (u)--(b1);
			\draw (a1)--(c1);
			\draw (b1)--(d1);
			\draw (c1)--(u)--(d1);
			\draw [line width = 1.2pt] (b1)--(c1)--(d1)--(b1);
			\end{scope}
			\begin{scope}[shift={(14,0)},scale=0.9]
			\node[vertex] (a1) at (3,1) {};
			\node[vertex] (b1) at (5.2,1) {};
			\node[vertex] (c1) at (5.2,-1.2) {};
			\node[vertex] (d1) at (3,-1.2) {};
			\node[vertex] (e1) at (7.4,1){};
			\node[vertex] (f1) at (7.4,-1.2){};
			\draw (a1)--(b1)--(c1)--(d1)--(a1)--(c1)--(e1)--(f1)--(b1)--(e1);
			\draw (f1)--(c1);
			\draw (d1)--(b1);
			\draw [line width = 1.2pt] (b1)--(c1);
			\end{scope}
			
			\end{tikzpicture}
			\caption{(left): $K_5^-$'s are edge-disjoint; (right): $K_4$'s are edge-disjoint. \label{fig:k5-k4edgedisjoint}}
		\end{figure}

		\item {\bf Step 3:} Cleaning $K_4$'s.
		
		First, we claim that any two copies of $K_4$ in $G$ are edge-disjoint. If not, then they can only intersect in one edge, or three edges. If they intersect in one edge, we find a $\widehat{P}_4$, and otherwise we get a $K_5^-$ in $G$ (see Figure \ref{fig:k5-k4edgedisjoint} (right); the intersecting edges are illustrated in bold). Hence, all $K_4$'s in $G$ are edge-disjoint.
		
		Consequently, if there are more than $\epsilon n^2$ copies of $K_4$ in $G$, taking one triangle from each copy gives us $\epsilon n^2$ edge-disjoint triangles, contradicting (\ref{eq:notManyEdgeDisjointTriangles}) again.
		
		Now, observe that for every $K_4$ in $G$ with vertices $\{a,b,c,d\}$, either the edge $ab$ or the edge $bc$ has codegree exactly $2$. Otherwise, let $x,y\in V(G)$ be such that $xab$ and $ybc$ are triangles in $G$. If $x=y$, then $xabcd$ is a $K_5^-$, and otherwise $xadcy$ is a $P_4$ in the neighborhood of $b$. Hence, whenever $G$ contains a $K_4$, we can remove an edge of codegree $2$ from it. Then $G$ loses at most $2\epsilon n^2$ triangles. Thus, we assume that $t(G)\ge \frac{n^2}{8}+5\epsilon n^2$, and that $G$ is $\{\widehat{P}_4, K_4 \}$-free.
		
		\item {\bf Step 4:} Cleaning $K_{2,2,2}$'s. 
		
		By assumption, $G$ contains no copy of $\widehat{P}_4$ and $K_4$. Fix a $K_{2,2,2}$ of $G$ with vertices $c_1,c_2$ in the center and $a_1,a_2,a_3,a_4$ forming the outer $C_4$. Let $X=V(G)\setminus\{c_1,c_2,a_1,a_2,a_3,a_4\}$. Denote $C_i = N_G(c_i)\cap X$ for $i=1,2$, and $A_i = N_G(a_i)\cap X$ for $i=1,\ldots,4$. Since $G$ is $\widehat{P}_4$-free, we deduce that $A_i\cap A_{i+1}=\varnothing$ and $A_i\cap C_j=\varnothing$ for every $i,j$ (here we denote $A_5:=A_1$). This is shown in Figure \ref{fig:w42block}, by assuming $x\in A_1\cap A_2$ and then $x\in A_1\cap C_2$, and finding copies of $\widehat{P}_4$ in either case. This implies that the $K_{2,2,2}$'s are themselves triangle blocks of $G$, hence they are mutually edge-disjoint.
		
		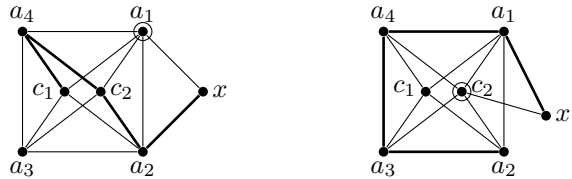
\begin{figure}[h]
			\centering
			\begin{tikzpicture}[scale = 0.8]
			\tikzstyle{vertex}=[circle,fill=black,minimum size=2pt,inner sep=1.3pt]
			\node[vertex] (a3) at (-1,1) {};
			\node[vertex] (b3) at (1,1) {};
			\node[vertex] (c3) at (1,-1) {};
			\node[vertex] (d3) at (-1,-1) {};
			\node[vertex] (u3) at (-0.3,0) {};
			\node[vertex] (v3) at (0.3,0) {};
			\node[vertex] (x3) at (2,0) {};
			\draw (a3)--(b3)--(c3)--(d3)--(a3);
			\draw (a3)--(u3)--(b3)--(v3)--(c3);
			\draw (c3)--(u3)--(d3)--(v3)--(a3);
			\draw (b3)--(x3)--(c3);
			\draw (b3) node[above] {$a_1$};
			\draw (c3) node[below] {$a_2$};
			\draw (d3) node[below] {$a_3$};
			\draw (a3) node[above] {$a_4$};
			\draw (u3) node[left] {$c_1$};
			\draw (v3) node[right] {$c_2$};						
			\draw (x3) node[right] {$x$};
			\draw (b3) circle (0.15cm);
			\draw [line width = 1pt] (x3)--(c3)--(v3)--(a3)--(u3);
			
			\begin{scope}[shift={(6,0)}]
			\node[vertex] (a3) at (-1,1) {};
			\node[vertex] (b3) at (1,1) {};
			\node[vertex] (c3) at (1,-1) {};
			\node[vertex] (d3) at (-1,-1) {};
			\node[vertex] (u3) at (-0.3,0) {};
			\node[vertex] (v3) at (0.3,0) {};
			\node[vertex] (x3) at (1.7,-0.4) {};
			\draw (a3)--(b3)--(c3)--(d3)--(a3);
			\draw (a3)--(u3)--(b3)--(v3)--(c3);
			\draw (c3)--(u3)--(d3)--(v3)--(a3);
			\draw (b3)--(x3)--(v3);
			\draw (b3) node[above] {$a_1$};
			\draw (c3) node[below] {$a_2$};
			\draw (d3) node[below] {$a_3$};
			\draw (a3) node[above] {$a_4$};
			\draw (-0.3,0.05) node[left] {$c_1$};
			\draw (0.3,0.05) node[right] {$c_2$};						
			\draw (x3) node[right] {$x$};
			\draw (v3) circle (0.15);
			\draw [line width = 1pt] (x3)--(b3)--(a3)--(d3)--(c3);
			\end{scope}
			\end{tikzpicture}
			\caption{\label{fig:w42block}$K_{2,2,2}$ is a block by itself.}
		\end{figure}
		
		So, if $G$ has at least $\epsilon n^2$ copies of $K_{2,2,2}$, then by taking one triangle from each $K_{2,2,2}$ we  obtain at least $\epsilon n^2$  edge-disjoint triangles in $G$, contradicting (\ref{eq:notManyEdgeDisjointTriangles}).
		
		Deleting one edge from each $K_{2,2,2}$, we lose at most $2\epsilon n^2$ triangles from $G$. Thus, we may assume that $t(G) \ge \frac{n^2}{8}+3\epsilon n^2$, and that $G$ is $\{\widehat{P}_4, K_4, K_{2,2,2}\}$-free.
		
		\item {\bf Step 5:} Cleaning $Q_{3,2}$'s.
		
		Suppose $G$ contains a $Q_{3,2}$ with the $P_3$ given by $a_1c_1c_2a_3$ and the outer $C_4$ being $a_1a_2a_3a_4$. Then, if $a_1c_2$ or $a_3c_1$ or $a_2a_4$ is an edge, we get a $K_4$ in $G$, and if $a_1a_3$ is an edge, then the $4$-cycle $a_1c_1c_2a_3$ along with vertices $a_2,a_4$ create a $K_{2,2,2}$ in $G$. Hence, every copy of $Q_{3,2}$ in $G$ has to be induced.
		
		Suppose $X=V(G)\setminus \{a_1,a_2,a_3,a_4,c_1,c_2\}$, and let $C_i = N_G(c_i)\cap X$ for $i=1,2$ and $A_i = N_G(a_i)\cap X$ for $i=1,\ldots,4$. Since $G$ is $\widehat{P}_4$-free, we deduce that $A_i\cap A_{i+1}=\varnothing$ and $A_i\cap C_j=\varnothing$ for every $i,j$ (here we denote $A_5:=A_1$), and $C_1\cap C_2=\varnothing$. We illustrate this in Figure \ref{fig:q32block}, similar to before. Hence, the $Q_{3,2}$'s of $G$ are themselves triangle blocks in $G$.
		
		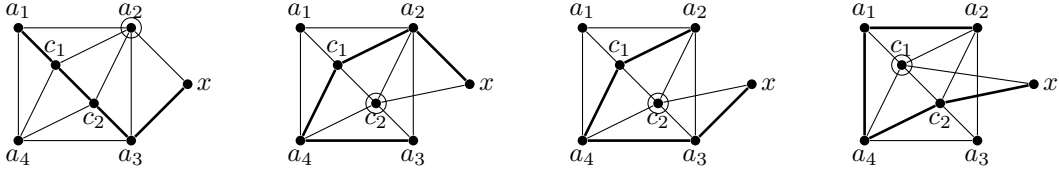
\begin{figure}[h]
			\centering
			\begin{tikzpicture}[scale=1.5]
			\tikzstyle{vertex}=[circle,fill=black,minimum size=2pt,inner sep=1.3pt]
			\node[vertex](a1)at(0,1){};
			\node[vertex](a2)at(1,1){};
			\node[vertex](a3)at(1,0){};
			\node[vertex](a4)at(0,0){};
			\node[vertex](c1)at($(a1)!0.33!(a3)$){};
			\node[vertex](c2)at($(a1)!0.67!(a3)$){};
			\node[vertex](x)at($(a2)!0.5!(a3) + (0.5,0)$){};
			\draw(a2)--(c1)--(a4)--(c2)--(a2)--(a1)--(a4)--(a3)--(a2);
			\draw(a1)--(a3);
			\draw(a2)circle(0.09);
			\draw[line width=1pt](a1)--(c1)--(c2)--(a3)--(x);
			\draw (x)--(a2);
			\draw(a1)node[above]{$a_1$};
			\draw(a2)node[above]{$a_2$};
			\draw(a3)node[below]{$a_3$};
			\draw(a4)node[below]{$a_4$};
			\draw(c1)node[above]{$c_1$};
			\draw(c2)node[below]{$c_2$};
			\draw(x)node[right]{$x$};
			
			\begin{scope}[shift={(2.5,0)}]
			\node[vertex](a1)at(0,1){};
			\node[vertex](a2)at(1,1){};
			\node[vertex](a3)at(1,0){};
			\node[vertex](a4)at(0,0){};
			\node[vertex](c1)at($(a1)!0.33!(a3)$){};
			\node[vertex](c2)at($(a1)!0.67!(a3)$){};
			\node[vertex](x)at($(a2)!0.5!(a3) + (0.5,0)$){};
			\draw(a2)--(c1)--(a4)--(c2)--(a2)--(a1)--(a4)--(a3)--(a2);
			\draw(a1)--(a3);
			\draw(c2)circle(0.09);
			\draw[line width=1pt](x)--(a2)--(c1)--(a4)--(a3);
			\draw(x)--(c2);
			\draw(a1)node[above]{$a_1$};
			\draw(a2)node[above]{$a_2$};
			\draw(a3)node[below]{$a_3$};
			\draw(a4)node[below]{$a_4$};
			\draw(c1)node[above]{$c_1$};
			\draw(c2)node[below]{$c_2$};
			\draw(x)node[right]{$x$};
			\end{scope}

			\begin{scope}[shift={(5,0)}]
			\node[vertex](a1)at(0,1){};
			\node[vertex](a2)at(1,1){};
			\node[vertex](a3)at(1,0){};
			\node[vertex](a4)at(0,0){};
			\node[vertex](c1)at($(a1)!0.33!(a3)$){};
			\node[vertex](c2)at($(a1)!0.67!(a3)$){};
			\node[vertex](x)at($(a2)!0.5!(a3) + (0.5,0)$){};
			\draw(a2)--(c1)--(a4)--(c2)--(a2)--(a1)--(a4)--(a3)--(a2);
			\draw(a1)--(a3);
			\draw(c2)circle(0.09);
			\draw[line width=1pt](x)--(a3)--(a4)--(c1)--(a2);
			\draw(c2)--(x);
			\draw(a1)node[above]{$a_1$};
			\draw(a2)node[above]{$a_2$};
			\draw(a3)node[below]{$a_3$};
			\draw(a4)node[below]{$a_4$};
			\draw(c1)node[above]{$c_1$};
			\draw(c2)node[below]{$c_2$};
			\draw(x)node[right]{$x$};
			\end{scope}
			
			\begin{scope}[shift={(7.5,0)}]
			\node[vertex](a1)at(0,1){};
			\node[vertex](a2)at(1,1){};
			\node[vertex](a3)at(1,0){};
			\node[vertex](a4)at(0,0){};
			\node[vertex](c1)at($(a1)!0.33!(a3)$){};
			\node[vertex](c2)at($(a1)!0.67!(a3)$){};
			\node[vertex](x)at($(a2)!0.5!(a3) + (0.5,0)$){};
			\draw(a2)--(c1)--(a4)--(c2)--(a2)--(a1)--(a4)--(a3)--(a2);
			\draw(a1)--(a3);
			\draw(c1)circle(0.09);
			\draw[line width=1pt](x)--(c2)--(a4)--(a1)--(a2);
			\draw(x)--(c1);
			\draw(a1)node[above]{$a_1$};
			\draw(a2)node[above]{$a_2$};
			\draw(a3)node[below]{$a_3$};
			\draw(a4)node[below]{$a_4$};
			\draw(c1)node[above]{$c_1$};
			\draw(c2)node[below]{$c_2$};
			\draw(x)node[right]{$x$};
			\end{scope}
		\end{tikzpicture}
		\caption{$Q_{3,2}$ is a block by itself.\label{fig:q32block}}
		\end{figure}
		
		Consequently, if $G$ has more than $\epsilon n^2$ copies of $Q_{3,2}$, then taking one triangle from each $Q_{3,2}$ we  obtain at least $\epsilon n^2$  edge-disjoint triangles in $G$, again contradicting (\ref{eq:notManyEdgeDisjointTriangles}). We delete an outer edge from each copy of $Q_{3,2}$, losing at most $\epsilon n^2$ triangles of $G$. Hence, we can assume that $t(G)\ge \frac{n^2}8+2\epsilon n^2$, and that $G$ is $\{\widehat P_4, K_4, K_{2,2,2},Q_{3,2}\}$-free.
		
		\item {\bf Step 6:} Cleaning $K_{1,2,2}$'s.
		
		We proceed similarly as before. First, for a $K_{1,2,2}$ with center $c$ and outer cycle $a_1a_2a_3a_4$, we call the edges $ca_i$ ``central edges". If $G$ contains such a $K_{1,2,2}$, then one cannot have an edge $a_1a_3$ or $a_2a_4$ since these give rise to $K_4$'s through $c$. Hence, the $K_{1,2,2}$'s in $G$ are induced. Plus, none of the edges $a_ic$ lie in a new triangle since it leads to a $\widehat{P}_4$: they all have codegree $2$. We now do a case analysis to see that the $K_{1,2,2}$'s in $G$ are edge-disjoint. Let $A,B\in \binom {V(G)}5$ be such that $G[A]$ and $G[B]$ are two $K_{1,2,2}$'s which are not edge-disjoint. Then $2\le |A\cap B|\le 4$. Let the central vertices of $G[A]$ and $G[B]$ be $u$ and $v$, respectively. Since each central edge of $G[A]$ and $G[B]$ has codegree $2$, $u\neq v$.
	
		Suppose $|A\cap B|=2$. If $u\in A\cap B$, then the edge through $u$ with its other endpoint in $A\cap B$ must have codegree at least $3$. Thus, the central vertices of $G[A]$ and $G[B]$ must lie outside $A\cap B$, leading us to the first configuration in Figure \ref{fig:wheelblocks}. But this configuration admits a $P_4$ in the neighborhood of either vertex of $A\cap B$, a contradiction. We illustrate $G[A\cap B]$ in boldface.
		
		\begin{figure}[h]
		\centering
			\begin{tikzpicture}[scale = 0.5]
			\tikzstyle{vertex}=[circle,fill=black,minimum size=2pt,inner sep=1.3pt]
			\begin{scope}[scale=1.2]
			\node[vertex] (a1) at (-1,1) {};
			\node[vertex] (a2) at (1,1) {};
			\node[vertex] (a3) at (1,-1) {};
			\node[vertex] (a4) at (-1,-1) {};
			\node[vertex] (c) at (0,0) {};
			\node[vertex] (b1) at (3,1){};
			\node[vertex] (b2) at (3,-1){};
			\node[vertex] (d) at (2,0) {};
			\draw (a2)circle(0.25);
			\draw (a1)--(a2)--(a3)--(a4)--(a1)--(c)--(a2);
			\draw (a3)--(c)--(a4);
			\draw (a3)--(d)--(a2);
			\draw (b1)--(d)--(b2);
			\draw (a2)--(b1)--(b2)--(a3);
			\draw (c)node[left]{$u$};
			\draw (d)node[right]{$v$};
			\draw [line width = 1.2pt] (a2)--(a3);
			\end{scope}
			\begin{scope}[shift={(8.5,0)}]
			\node[vertex] (a1) at (-2,0) {};
			\node[vertex] (c) at (-1,0) {};
			\node[vertex] (a3) at (0,0) {};
			\node[vertex] (a4) at (0,-2) {};
			\node[vertex] (a2) at (0,2) {};
			\node[vertex] (c1) at (1,0) {};
			\node[vertex] (b1) at (2,0) {};
			\draw (a2)circle(0.27);
			\draw (a1)--(a2)--(b1)--(a4)--(a1)--(c)--(a2)--(c1)--(a3)--(c)--(a4)--(c1)--(b1);
			\draw [line width = 1.2pt] (a2)--(a3)--(a4);
			\draw (c)node[below]{$u\ $};
			\draw (c1)node[below]{$\ v$};
			\end{scope}
			\begin{scope}[shift={(16,0)},scale = 1.2]
			\node[vertex] (a1) at (-1,1) {};
			\node[vertex] (a2) at (1,1) {};
			\node[vertex] (a3) at (1,-1) {};
			\node[vertex] (a4) at (-1,-1) {};
			\node[vertex] (c) at (-2,0) {};
			\node[vertex] (d) at (2,0) {};
			\draw (c)node[left]{$u$};
			\draw (d)node[right]{$v$};
			\draw (a1)--(c)--(a2)--(d)--(a3)--(c)--(a4)--(d)--(a1);
			\draw [line width = 1.2pt] (a1)--(a2)--(a3)--(a4)--(a1);
			\end{scope}
			\begin{scope}[shift={(23,0)}]
			\node[vertex] (a) at (-2,0) {};
			\node[vertex] (c) at (2,0) {};
			\node[vertex] (d) at (0,-2) {};
			\node[vertex] (b) at (0,2) {};
			\node[vertex] (x) at ($(a)!0.33!(c)$) {};
			\node[vertex] (y) at ($(a)!0.67!(c)$) {};
			\draw (x)node[below]{$u\ \ $};
			\draw (y)node[below]{$\ \ v$};
			\draw (a)--(b)--(c)--(d)--(a)--(x)--(b)--(y)--(x)--(d)--(y);
			\draw (y)--(c);
			\draw [line width = 1.2pt] (x)--(b)--(y)--(d)--(x)--(y);
			\end{scope}
			\end{tikzpicture}
		\caption{\label{fig:wheelblocks}The different ways two induced $K_{1,2,2}$'s can intersect.}
		\end{figure}
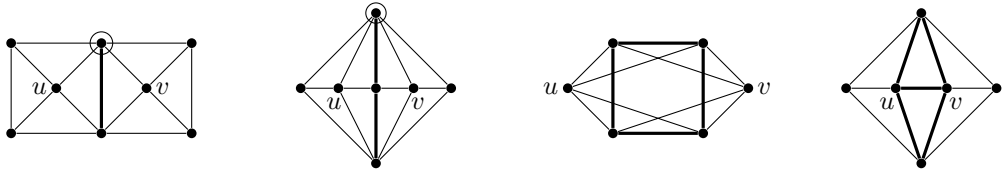
		
		Next, suppose $|A\cap B|=3$.  If $u\in A\cap B$ but $v\not\in A\cap B$, then one of the central edges of $G[A]$ contains an external triangle through $v$. If $u,v\in A\cap B$, then as $u$ is part of the outer $C_4$ of $G[B]$, $uv$ contains a triangle from $G[B]$ which is not contained in $G[A]$. Thus the only possibility for $|A\cap B|=3$ is for $u$ and $v$ to be both outside $A\cap B$. This gives rise to the second configuration in Figure \ref{fig:wheelblocks}, which contains $\widehat P_4$ in the neighborhood of one of the vertices of $A\cap B$.
		
		Finally, if $|A\cap B|=4$ and $u\in A\cap B$ but $v\not\in A\cap B$, then any edge of $G[A\cap B]$ through $u$ has codegree at least $3$. If both $u$ and $v$ lie outside $A\cap B$, we obtain a $K_{2,2,2}$, which is the third configuration in Figure \ref{fig:wheelblocks}. Hence, $u$ and $v$ must both lie inside $A\cap B$. Since $u$ and $v$ both must be adjacent to all other vertices of $A\cap B$, $G[A\cup B]$ form a $Q_{3,2}$, the fourth configuration in Figure \ref{fig:wheelblocks}.
		
		Therefore, if two $K_{1,2,2}$'s are not edge-disjoint, they must intersect each other in one of the ways depicted in Figure \ref{fig:wheelblocks}, and we either find a $\widehat P_4$, $K_{2,2,2}$ or $Q_{3,2}$ inside $G$  for each of these intersecting patterns. Thus, all $K_{1,2,2}$'s of $G$ are edge-disjoint. Consequently, if $G$ has $\epsilon n^2$ copies of $K_{1,2,2}$, they are all edge-disjoint, and give us at least $\epsilon n^2$ edge-disjoint triangles, again contradicting (\ref{eq:notManyEdgeDisjointTriangles}).
		
		For each $K_{1,2,2}$ in $G$ with central vertex $x$ and outer cycle $abcd$, we observe that either $ab$ or $bc$ has codegree $1$. Otherwise, suppose $y,z\in V(G)$ are such that $yab$ and $zbc$ form triangles in $G$. If $y=z$, this creates a $Q_{3,2}$ in $G$. Otherwise, $yaxcz$ is a $P_4$ in the neighborhood of $b$. So, every $K_{1,2,2}$ has an outer edge of codegree $1$.	By deleting one such edge of codegree $1$ from each copy of $K_{1,2,2}$, we remove at most $\epsilon n^2$ triangles from $G$.
		Therefore, we may assume that $G$ is $\{\widehat{P}_4, K_4, K_{1,2,2}\}$-free, and 
		\[
		t(G)\ge \frac{n^2}8 + \epsilon n^2.
		\]
	\end{itemize}
	\medskip

		Let us now analyze the structure of $G$. We will prove using induction on $t(H)$, that for any subgraph $H\subseteq G$,
		\begin{equation}
		\label{eq:inductionP4}
		t(H)\le\frac{e(H)}2.
		\end{equation}
		
		When $t(H)=1$, $e(H)\ge 3$, proving the base case. Now suppose $t(H)>1$ for some $H\subseteq G$, and that (\ref{eq:inductionP4}) holds for all subgraphs $H'$ with $t(H')<t(H)$. Assume without loss of generality that every edge of $H$ lies in at least one triangle. Call an edge of $H$ \emph{light} if it is contained in a unique triangle from $H$. Call edges that are not light, \emph{heavy}. We observe that if $H$ contains a triangle with two light edges, then deleting them from $H$ leads to a graph $H'\subsetneq H$ with $t(H')=t(H)-1$ and $e(H')=e(H)-2$. Using the induction hypothesis on $H'$, $t(H')\le e(H')/2$, implying $t(H)\le e(H)/2$. Hence, we may further assume that $H$ contains no triangle with two light edges.
		
		\begin{lemma}
			\label{lem:blockneighborhood}
			Suppose $H$ contains two triangles $xuv$ and $yuv$ intersecting in the edge $uv$. Then either: (a) $xu,yv$ are light and $xv,yu$ are heavy or: (b) $xu,yv$ are heavy and $xv,yu$ are light.
		\end{lemma}
		\begin{figure}[h]
			\centering
			\begin{tikzpicture}[scale=0.7]
			\tikzstyle{vertex}=[circle,fill=black,minimum size=2pt,inner sep=1.3pt]
			\node[vertex] (y0) at (1,1){};
			\node[vertex] (x0) at (0,0){};
			\node[vertex] (z0) at (2,0){};
			\node[vertex] (v0) at (1,-1){};
			\draw (x0)--(y0)--(z0)--(v0)--(x0);
			\draw (v0)--(y0);
			\draw [dashed] (x0)--(z0);
			\draw (x0) node [left] {$x$};
			\draw (y0) node [above] {$u$};
			\draw (z0) node [right] {$y$};
			\draw (v0) node [below] {$v$};
			
			\begin{scope}[shift={(4,0)}]
			\node[vertex] (y1) at (5,1){};
			\node[vertex] (x1) at (4,0){};
			\node[vertex] (z1) at (6,0){};
			\node[vertex] (v1) at (5,-1){};
			\node[vertex] (u1) at (6,1){};
			\draw (x1)--(y1)--(z1)--(v1);
			\draw (v1)--(y1);
			\draw (v1)--(x1);
			\draw (u1)--(x1);
			\draw (u1)--(y1);
			\draw (u1)--(z1);
			\draw (x1) node [left] {$x$};
			\draw (y1) node [above] {$u$};
			\draw (z1) node [right] {$y$};
			\draw (v1) node [below] {$v$};
			\draw (u1) node [right] {$z_1=z_2$};
			\end{scope}
			
			\begin{scope}[shift={(-6,0)}]
			\node[vertex] (y2) at (11,1){};
			\node[vertex] (x2) at (10,0){};
			\node[vertex] (z2) at (12,0){};
			\node[vertex] (v2) at (11,-1){};
			\node[vertex] (u22) at (12,1){};
			\node[vertex] (u21) at (10,1){};
			\draw (x2)--(y2)--(z2)--(v2);
			\draw (v2)--(y2);
			\draw (v2)--(x2);
			\draw (u21)--(x2);
			\draw (u21)--(y2);
			\draw (u22)--(y2);
			\draw (u22)--(z2);
			\draw (x2) node [left] {$x$};
			\draw (y2) node [above] {$u$};
			\draw (z2) node [right] {$y$};
			\draw (v2) node [below] {$v$};
			\draw (u21) node [left] {$z_1$};
			\draw (u22) node [right] {$z_2$};
			\draw (y2) circle (0.15);
			\draw [line width = 1pt] (u22)--(z2)--(v2)--(x2)--(u21);
			\end{scope}
			\end{tikzpicture}
			\caption{\label{fig:blockneighborhood}$xu$ and $yu$ cannot both be heavy.}
		\end{figure}
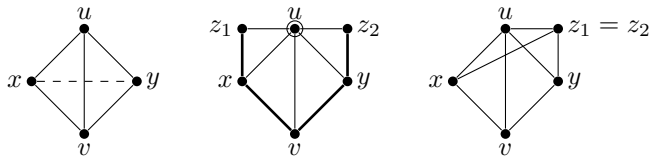
		\noindent{\it Proof of Lemma \ref{lem:blockneighborhood}.}
		Suppose that both $xu$ and $yu$ were heavy (Figure \ref{fig:blockneighborhood}). If $x$ and $y$ were adjacent, this would create a $K_4$ which is forbidden. So, there exist $z_1,z_2\in V(H)$ such that $z_1xu$ and $z_2yu$ form $K_3$'s in $H$. If $z_1\neq z_2$, $N_H(u)$ contains a $P_4$, which is forbidden. Otherwise $z_1=z_2$, and this produces a $K_{1,2,2}$ centered at $u$, a contradiction.
		
		Hence one of $xu$ and $yu$ is light. Similarly, one of $xv$ and $yv$ is light. If $xu$ is light, then $xv$ and $yu$ are heavy, implying that $yv$ is light, and \emph{(a)} holds. Similarly, if $xu$ is heavy, then \emph{(b)} holds.
		
\begin{figure}[h]
			\centering
			\begin{tikzpicture}[scale=0.7]
			\tikzstyle{vertex}=[circle,fill=black,minimum size=2pt,inner sep=1.3pt]
			\node[vertex] (x1) at (0,2){};
			\node[vertex] (x2) at (2,2){};
			\node[vertex] (y1) at (1,0){};
			\node[vertex] (y2) at (3,0){};
			\node[vertex] (t) at (4,2){};
			\draw (x1) node [above] {$x$};
			\draw (x2) node [above] {$u$};
			\draw (y1) node [below] {$v$};
			\draw (y2) node [below] {$y$};
			\draw (t) node [above] {$z$};
			\draw [dashed](x1)--(t);
			\draw [dashed](y2)--(y1)--(t);
			\draw (x1)--(y1)--(x2)--(y2);
		\end{tikzpicture}
		\caption{\label{fig:codegree<3}Codegree of $uv\in E(H)$ is at most $2$.}
		\end{figure}
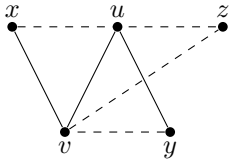

		\medskip
		We shall now use Lemma \ref{lem:blockneighborhood} and the fact that every triangle in $H$ has two heavy and one light edge, to analyze the structure of $H$. First, observe that $H$ cannot have any edge of codegree more than $2$. This is because if we have an edge $uv$ which lies in three triangles $xuv$, $yuv$, $zuv$, then by Lemma \ref{lem:blockneighborhood}, either $xu, yv$ are light or $xv, yu$ are light. Suppose without loss of generality that $xu$ and $yv$ are light, as in Figure \ref{fig:codegree<3}. Then, by applying Lemma \ref{lem:blockneighborhood} on the pairs $\{xuv,zuv\}$ and $\{yuv,zuv\}$ respectively, the edges $zv$ and $zu$ must be light. However, this contradicts the assumption of $H$ containing no triangle with two light edges.
			
		Now, let $\ell(H)$ denote the number of light edges of $H$ and $h(H)$ the number of heavy edges of $H$.	Since every edge of $H$ can have codegree $1$ or $2$, and every triangle contains one light and two heavy edges, a double-counting argument gives,
		\[
		\ell(H)+2h(H) = 3t(H).
		\]
		On the other hand, every light edge of $H$ lies in a unique triangle, and every triangle contains a unique light edge. This implies $t(H)=\ell(H)$. Therefore,
		\[
		4t(H) = 2\ell(H) + 2 h(H) = 2e(H),
		\]
		implying $t(H)=\frac{e(H)}2$. This finishes the induction step, completing the proof of (\ref{eq:inductionP4}).
		
		\medskip
		Taking $H=G$ in (\ref{eq:inductionP4}), we obtain $t(G)\le \frac{e(G)}2$. By assumption, $t(G)\ge \frac{n^2}{8}+\epsilon n^2$. So, by Lemma \ref{lem:nordhausStewart},
		\[
		\frac{e(G)}{2}\ge t(G)\ge \frac{e(G)}{3n}\cdot (4e(G)-n^2),
		\]
		leading to $e(G)\le \frac{n^2}{4}+\frac{3n}{8}$. This gives $t(G)\le \frac{e(G)}{2}\le \frac{n^2}{8}+\frac{3n}{16}$, which contradicts the assumption of $t(G)\ge \frac{n^2}8 + \epsilon n^2$ for sufficiently large $n$.		This concludes the proof of Theorem \ref{thm:p3p4p5hat} for $k=4$.
\hfill{$\Box$}

\subsubsection{Proof of Theorem \ref{thm:p3p4p5hat} for $\mathbf{k=5}$.}
	Our proof of Theorem \ref{thm:p3p4p5hat} for $k=5$ follows exactly the same structure as that for $k=3$ and $k=4$, with more technical details.
	We shall prove, using induction on $n$, that if $G$ is $\widehat P_5$-free, then 
	\begin{equation}
	\label{eq:inductionP5onN}
	t(G)\le \frac{n^2}{4} + 5n.
	\end{equation}

	The base case $n=3$ is clearly true as $t(G)\le 1$. Assume that (\ref{eq:inductionP5onN}) holds for all graphs $G$ on less than $n$ vertices, and let us prove that it also holds for $G$. The first step is to remove all copies of $K_6$ and $K_6^-$ from $G$ via the induction hypothesis.
	
	Suppose $G$ has a copy of $K_6$ with vertices $a_1,\ldots, a_6$. Then this is a block by itself, since if there is a vertex $x\neq a_i$ such that $xa_1a_2$ is a triangle, then $N_G(a_1)$ contains the $5$-path $a_6a_5a_4a_3a_2x$. For $1\le i\le 6$, let $X_i=N_G(a_i)\setminus\{a_1,\ldots,a_6\}$. Then $X_i\cap X_j=\varnothing$ for every $i\neq j$. Since $\sum_{i=1}^6|X_i|\le n-6$, there is a vertex $a_i$ for which $|X_i|\le \frac{n-6}{6}$. By Lemma~\ref{lem:erdosGallaiTheta},
	\[
	e(X_i)\le \frac{5-1}{2}\cdot |X_i|\le \frac{n-6}{3}.
	\]
	Hence, by (\ref{eq:inductionP5onN}) on $G'=G-\{a_i\}$, we get $t(G')\le \frac{(n-1)^2}{4}+5(n-1)$. Therefore,
	\[
	t(G)\le t(G') + \frac{n-6}{3}+5\le \frac{(n-1)^2}{4}+5(n-1)+\frac{n-6}{3}+5 < \frac{n^2}{4}+5n,
	\]
	completing the induction step for $G$. Hence, we may assume that $G$ is $K_6$-free.
	
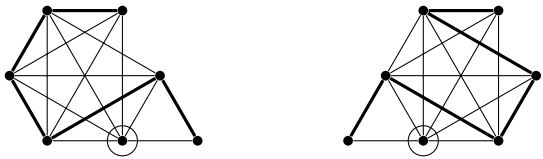
\begin{figure}[h]
		\centering
		\begin{tikzpicture}
		\tikzstyle{vertex}=[circle,fill=black,minimum size=2pt,inner sep=1.3pt]
		\node[vertex](a1)at(1,0){};
		\node[vertex](a2)at(60:1){};
		\node[vertex](a3)at(120:1){};
		\node[vertex](a4)at(180:1){};
		\node[vertex](a5)at(240:1){};
		\node[vertex](a6)at(300:1){};
		\draw(a1)--(a3)--(a4)--(a5)--(a6)--(a2)--(a3)--(a5)--(a1)--(a4)--(a2)--(a5);
		\draw(a1)--(a6)--(a3);
		\draw(a6)--(a4);
		\node[vertex](x)at(330:1.732){};
		\draw(a1)--(x)--(a6);
		\draw(a6)circle(0.2);
		\draw[line width=1.2pt] (x)--(a1)--(a5)--(a4)--(a3)--(a2);
		
		\begin{scope}[shift={(5,0)}]
		\node[vertex](a1)at(1,0){};
		\node[vertex](a2)at(60:1){};
		\node[vertex](a3)at(120:1){};
		\node[vertex](a4)at(180:1){};
		\node[vertex](a5)at(240:1){};
		\node[vertex](a6)at(300:1){};
		\draw(a1)--(a3)--(a4)--(a5)--(a6)--(a2)--(a3)--(a5)--(a1)--(a4)--(a2)--(a5);
		\draw(a1)--(a6)--(a3);
		\draw(a6)--(a4);
		\node[vertex](x)at(210:1.732){};
		\draw(a5)--(x)--(a4);
		\draw(a5)circle(0.2);
		\draw[line width=1.2pt] (x)--(a4)--(a6)--(a1)--(a3)--(a2);
		\end{scope}
		\end{tikzpicture}
		\caption{$K_6^-$ is a block by itself.\label{fig:k6minusBlock}}
	\end{figure}

	Now, if $G$ has a copy of $K_6^-$ on vertices $a_1,\ldots,a_6$, it has to be induced. We verify in Figure~\ref{fig:k6minusBlock} that it is a block by finding a $\widehat P_5$ whenever any edge lies in an external triangle. Let $X_i=N_G(a_i)\setminus\{a_1,\ldots,a_6\}$. Following exactly the same argument as before, there exists a vertex $a_i$ for which $|X_i|\le \frac{n-6}{6}$. Therefore, applying Lemma~\ref{lem:erdosGallaiTheta} and letting $G'=G-\{a_i\}$,
	\[
	t(G)\le t(G') + \frac{n-6}{3} + 5 <\frac{n^2}{4}+5n,
	\]
	completing the induction step for $G$.
	
	Therefore, without loss of generality we can assume that $G$ is $\{K_6^-,\widehat P_5\}$-free. Consider $G$ to be a fixed $n$-vertex graph. We shall now prove using induction on $t(H)$, that for any subgraph $H\subseteq G$,
	\begin{equation}
	\label{eq:inductionP5}
	t(H)\le e(H).
	\end{equation}
	
	When $t(H)=1$, $e(H)\ge 3$ proves the base case. Now suppose $t(H)>1$ for some $H\subseteq G$, and that (\ref{eq:inductionP5}) holds for all subgraphs $H'$ of $G$ with $t(H')<t(H)$. If $H$ has an edge $e$ which lies in at most one triangle, using the induction hypothesis on $H'=H-\{e\}$ immediately proves (\ref{eq:inductionP5}) for $H$. Hence, we may assume that all edges of $H$ have codegree at least $2$. Call an edge of $H$ \emph{light} if it has codegree exactly $2$, otherwise call it \emph{heavy}.
	
	\begin{lemma}
		\label{lem:cleaningForP5}
		We may assume that $H$ does not contain $W_5$ and $K_{1,2,2}$ as subgraphs.
	\end{lemma}
	
	This lemma is proved by sequentially removing copies of the graphs illustrated in Figure~\ref{fig:cleaningListForP5} from $H$, and the proof can be found in the Appendix. We shall now assume that Lemma~\ref{lem:cleaningForP5} is true.
	
	\begin{figure}[h]
		\centering
		\begin{tikzpicture}[scale=0.7]
		\tikzstyle{vertex}=[circle,fill=black,minimum size=2pt,inner sep=1.3pt]
		\node[vertex](a1)at(1,0){};
		\node[vertex](a2)at(60:1){};
		\node[vertex](a3)at(120:1){};
		\node[vertex](a4)at(180:1){};
		\node[vertex](a5)at(240:1){};
		\node[vertex](a6)at(300:1){};
		\draw(a1)--(a3)--(a4)--(a5)--(a6)--(a2)--(a5)--(a1)--(a4)--(a2)--(a5);
		\draw(a1)--(a6)--(a3)--(a5);
		\draw(a6)--(a4);
		\draw(0,-1.5)node{$K_6^{-2,1}$};
		
		\begin{scope}[shift={(3,0)}]
		\node[vertex](a1)at(1,0){};
		\node[vertex](a2)at(60:1){};
		\node[vertex](a3)at(120:1){};
		\node[vertex](a4)at(180:1){};
		\node[vertex](a5)at(240:1){};
		\node[vertex](a6)at(300:1){};
		\draw(a1)--(a3)--(a5)--(a6)--(a2)--(a3)--(a5)--(a1)--(a4)--(a2)--(a5);
		\draw(a1)--(a6)--(a3);
		\draw(a6)--(a4)--(a5);
		\draw(0,-1.5)node{$K_6^{-2,2}$};
		\end{scope}
		
		\begin{scope}[shift={(6,0)}]
		\node[vertex](a1)at(1,0){};
		\node[vertex](a2)at(60:1){};
		\node[vertex](a3)at(120:1){};
		\node[vertex](a4)at(180:1){};
		\node[vertex](a5)at(240:1){};
		\node[vertex](a6)at(300:1){};
		\draw(a1)--(a3)--(a5)--(a6)--(a2)--(a5)--(a1)--(a4)--(a2)--(a5);
		\draw(a1)--(a6)--(a3)--(a5);
		\draw(a6)--(a4)--(a5);
		\draw(0,-1.5)node{$K_6^{-3,1}$};
		\end{scope}

		\begin{scope}[shift={(9,0)}]
		\node[vertex](a1)at(1,0){};
		\node[vertex](a2)at(60:1){};
		\node[vertex](a3)at(120:1){};
		\node[vertex](a4)at(180:1){};
		\node[vertex](a5)at(240:1){};
		\node[vertex](a6)at(300:1){};
		\draw(a1)--(a3)--(a5)--(a6)--(a2)--(a5)--(a1)--(a4)--(a2)--(a5);
		\draw(a3)--(a4)--(a6)--(a1);
		\draw(a3)--(a6);
		\draw(0,-1.5)node{$K_6^{-3,2}$};
		\end{scope}
		
		\begin{scope}[shift={(12,0)},rotate=144]
		\node[vertex](a1)at(1,0){};
		\node[vertex](a2)at(72:1){};
		\node[vertex](a3)at(144:1){};
		\node[vertex](a4)at(216:1){};
		\node[vertex](a5)at(288:1){};
		\draw(a2)--(a3)--(a4)--(a5)--(a1)--(a3)--(a5)--(a2)--(a4)--(a1);
		\end{scope}
		\draw(12,-1.5)node{$K_5^{-}$};
		
		\begin{scope}[shift={(15,0)},rotate=144]
		\node[vertex](a1)at(1,0){};
		\node[vertex](a2)at(72:1){};
		\node[vertex](a3)at(144:1){};
		\node[vertex](a4)at(216:1){};
		\node[vertex](a5)at(288:1){};
		\node[vertex](x) at(0,0){};
		\draw(a1)--(a2)--(a3)--(a4)--(a5)--(a1)--(x)--(a3)--(a5)--(x)--(a4);
		\draw(x)--(a2);
		\end{scope}
		\draw(15,-1.5)node{$W_5^{+}$};
		
		\begin{scope}[shift={(18,0)}]
		\node[vertex](a1)at(1,0){};
		\node[vertex](a2)at(72:1){};
		\node[vertex](a3)at(144:1){};
		\node[vertex](a4)at(216:1){};
		\node[vertex](a5)at(288:1){};
		\node[vertex](x) at(0,0){};
		\draw(a1)--(a2)--(a3)--(a4)--(a5)--(a1)--(x)--(a3)--(x)--(a4);
		\draw(a5)--(x)--(a2);
		\draw(0,-1.5)node{$W_5$};
		\end{scope}
		
		\begin{scope}[shift={(21,0)},rotate=45]
		\node[vertex](a1)at(1,0){};
		\node[vertex](a2)at(90:1){};
		\node[vertex](a3)at(180:1){};
		\node[vertex](a4)at(270:1){};
		\node[vertex](x) at(0,0){};
		\draw(a1)--(a2)--(a3)--(a4)--(x)--(a1)--(a4);
		\draw(a3)--(x)--(a2);
		\end{scope}
		\draw(21,-1.5)node{$W_4=K_{1,2,2}$};
		\end{tikzpicture}
		\caption{Graphs to be cleaned from $H$.\label{fig:cleaningListForP5}}
	\end{figure}
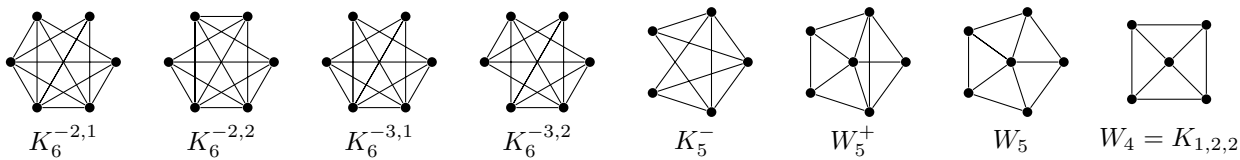

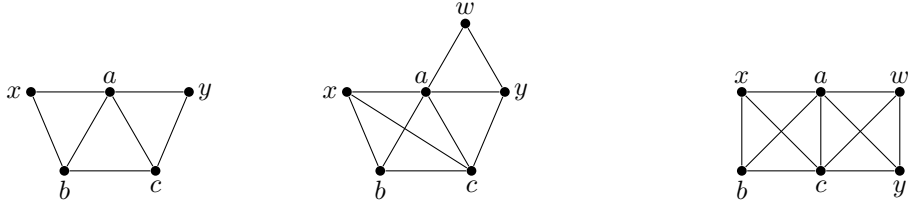
\begin{figure}[h]
	\centering
	\begin{tikzpicture}[scale=0.7]
	\tikzstyle{vertex}=[circle,fill=black,minimum size=2pt,inner sep=1.3pt]
	\node[vertex](a)at(0,1){};
	\node[vertex](b)at(210:1){};
	\node[vertex](c)at(-30:1){};
	\node[vertex](x)at(-1.5,1){};
	\node[vertex](y)at(1.5,1){};
	\draw (b)--(x)--(a)--(y)--(c)--(b)--(a)--(c);
	\draw (a)node[above]{$a$};
	\draw (b)node[below]{$b$};
	\draw (c)node[below]{$c$};
	\draw (x)node[left]{$x$};
	\draw (y)node[right]{$y$};
	
	\begin{scope}[shift={(6,0)}]
	\node[vertex](a)at(0,1){};
	\node[vertex](b)at(210:1){};
	\node[vertex](c)at(-30:1){};
	\node[vertex](x)at(-1.5,1){};
	\node[vertex](y)at(1.5,1){};
	\node[vertex](w)at(0.75,2.3){};
	\draw (b)--(x)--(a)--(y)--(c)--(b)--(a)--(c);
	\draw (x)--(c);
	\draw (a)--(w)--(y);
	\draw (a)node[above]{$a\ $};
	\draw (b)node[below]{$b$};
	\draw (c)node[below]{$c$};
	\draw (x)node[left]{$x$};
	\draw (y)node[right]{$y$};
	\draw (w)node[above]{$w$};
	\end{scope}
	
	\begin{scope}[shift={(12,0)},xscale=1.5]
	\node[vertex](a)at(1,1){};
	\node[vertex](b)at(0,-0.5){};
	\node[vertex](c)at(1,-0.5){};
	\node[vertex](x)at(0,1){};
	\node[vertex](y)at(2,-0.5){};
	\node[vertex](w)at(2,1){};
	\draw (b)--(x)--(a)--(y)--(c)--(b)--(a)--(c);
	\draw (x)--(c)--(w);
	\draw (a)--(w)--(y);
	\draw (a)node[above]{$a$};
	\draw (b)node[below]{$b$};
	\draw (c)node[below]{$c$};
	\draw (x)node[above]{$x$};
	\draw (y)node[below]{$y$};
	\draw (w)node[above]{$w$};
	\end{scope}
	\end{tikzpicture}
	\caption{$abc$, $abx$, $acy$ are triangles with $x\neq y$.\label{fig:abcxyTriangles}}
\end{figure}

	Suppose $H$ contains a triangle $abc$ such that $abx$ and $acy$ are triangles, with $x\neq y$, i.e $H$ contains a $\widehat P_3$. Refer to Figure~\ref{fig:abcxyTriangles}. Observe that both $ax$ and $ay$ must have codegree at least $2$. If $xy$ is an edge in $H$, we get a $W_4$. If $axz$ and $ayw$ are triangles for vertices $z$ and $w$ which are not $b$ or $c$, then there are two possibilities. Either $z\neq w$ in which case we get a $\widehat P_5$, or $z=w$, producing a $W_5$ in $H$. Therefore $\{z,w\}\cap \{b,c\}\neq\varnothing$. If $z=c$ and $w=b$, this gives us a $K_{1,2,2}$ centered at $a$. Hence we may assume $z=c$ and $w\neq b$. By assumption, $aw$ must have codegree at least $2$. Note that $wb$ or $wx$ cannot be edges, as they create $K_{1,2,2}$ or $W_5$ in $H$ centered around $a$, respectively. Further, we cannot have a new vertex $t$ for which $awt$ is a triangle, since this creates a $\widehat P_5$ centered at $a$. Thus, the only possibility is that $wc\in E(H)$.

	As $H$ is $\{K_{1,2,2},W_5\}$-free, $H[a,b,c,x,y,w]$ is induced. Further, if the edges $ax,ab,ay,aw$ or $cb,cx,cw,cy$ lie in an external triangle, we can find $\widehat{P}_5$'s centered around $a$ or $c$, respectively. Hence these $8$ edges all do not lie in external triangles, and have codegree exactly $2$. Deleting them from $H$, we obtain a graph $H'$ with $t(H')=t(H)-8$ and $e(H')=e(H)-8$, completing the proof of (\ref{eq:inductionP5}) for $H$.

	Hence, we may assume that $H$ does not contain any $\widehat P_3$. Now if $ab$ and $ac$ were heavy in any triangle $abc$, we would then find $x\neq y$ for which $abx$ and $acy$ are triangles in $H$. This leads us to the following crucial observation:
	\begin{equation}
	\label{eq:oneHeavyEdgePerTriangle}
	\mbox{Every triangle of }H\mbox{ has at most one heavy edge.}
	\end{equation}
	
	\begin{figure}[h]
		\centering
		\begin{tikzpicture}
			\tikzstyle{vertex}=[circle,fill=black,minimum size=2pt,inner sep=1.3pt]
			\node[vertex](a)at(90:1){};
			\node[vertex](b)at(210:1){};
			\node[vertex](c)at(-30:1){};
			\node[vertex](x)at(0,0){};
			\draw[dashed](b)--(a)--(c);
			\draw[dashed](x)--(a);
			\draw(x)--(b)--(c)--(x);
			\draw(a)node[above]{$a$};
			\draw(b)node[left]{$b$};
			\draw(c)node[right]{$c$};
			\draw(x)node[below]{$x$};
			\begin{scope}[shift={(4,0)}]
				\node[vertex](a)at(90:1){};
				\node[vertex](b)at(210:1){};
				\node[vertex](c)at(-30:1){};
				\node[vertex](x)at(0,0){};
				\draw[dashed](b)--(a)--(c);
				\draw(x)--(a);
				\draw[dashed](b)--(x)--(c);
				\draw(b)--(c);
				\draw(a)node[above]{$a$};
				\draw(b)node[left]{$b$};
				\draw(c)node[right]{$c$};
				\draw(x)node[below]{$x$};
			\end{scope}
		\end{tikzpicture}
		\caption{\label{fig:P5hatMainProof}Edge $xa$ can be light or heavy.}
	\end{figure}
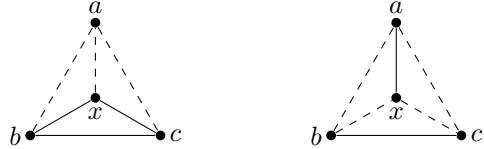

	Let us fix a triangle $abc$ in $H$. Let $ab$ and $ac$ be light. As they must have codegree $2$, there is a vertex $x$ for which $xa,xb,xc\in E(H)$, as in Figure~\ref{fig:P5hatMainProof}. If the edge $xa$ is light, we can then let $H'=H-\{ab,ax,ac\}$. Note that $t(H')=t(H)-3$ and $e(H')=e(H)-3$, finishing the proof of (\ref{eq:inductionP5}) for $H$. Finally, if $xa$ is heavy, then by (\ref{eq:oneHeavyEdgePerTriangle}), the edges $xb$ and $xc$ must be light. Let $H'=H-\{ab,ac,xb,xc\}$, then $t(H')=t(H)-4$ and $e(H')=e(H)-4$, completing the induction step of (\ref{eq:inductionP5}).

	Taking $H=G$ in (\ref{eq:inductionP5}), we obtain $t(G)\le e(G)$. Using Lemma~\ref{lem:nordhausStewart},
	\[
	e(G)\ge t(G)\ge \frac{e(G)}{3n}\cdot(4e(G)-n^2),
	\]
	implying $t(G)\le e(G)\le \frac{n^2}{4}+\frac{3n}{4}<\frac{n^2}{4}+5n$. This concludes the proof of (\ref{eq:inductionP5}) for $G$.
	
\section{Concluding Remarks}
The following generalization of Theorems \ref{thm:k1ab} and \ref{thm:w2k} was obtained by Abhishek Methuku during the review process.

\begin{thm}
	\label{thm:reviewerExtra}
	Suppose $H$ is a graph with $\ex(n,H)=O(n^\alpha)$ for some $1< \alpha < 2$. Then, $\ex(n, K_3, \widehat H) = o(n^{1+\alpha})$.
\end{thm}

\begin{proof}[Proof of Theorem \ref{thm:reviewerExtra}]
	Recall that $H$ is a graph such that $\ex(n,H)\le Cn^\alpha$, where $C>0, 1<\alpha<2$, and we want to show that $\ex(n,K_3,\widehat H)=o(n^{1+\alpha})$.
	First, let us fix a graph $G$ on $n$ vertices that is $\widehat H$-free, and any $\epsilon>0$.
	The number of triangles in $G$ is at most $n\cdot \ex(n,H) \le Cn^{1+\alpha} = o(n^3)$ many triangles.
	Thus, by Lemma \ref{lem:triangleRemoval}, every maximal collection of edge-disjoint triangles in $G$ has size at most $\epsilon n^2$.
	Hence, there is a set $F$ of edges of $G$ such that $|F|=3\epsilon n^2$ and every triangle of $G$ has an edge in $F$.
	For every vertex $v\in V(G)$, define $N_F(v) := \{x\in V(G): vx \in F\} \subseteq N_G(v)$, and $N_R(v) := N_G(v)\setminus N_F(v)$.
	Since every triangle of $G$ is incident to $N_F(v)$ for some $v\in V(G)$,
	\begin{equation}
	\label{eq:concl1}
	t(G)\le \sum_{v\in V(G)} e(N_F(v)) + \sum_{v\in V(G)} e(N_F(v), N_R(v)).
	\end{equation}
	Here $e(N_F(v), N_R(v))$ denotes the number of edges between vertex subsets $N_F(v)$ and $N_R(v)$ in $G$.
	Now we bound each term in (\ref{eq:concl1}).
	By assumption, as $N_G(v)$ is $H$-free, we have 
	\[
	e(N_F(v))\le \ex(|N_F(v)|,H)\le C|N_F(v)|^{\alpha}\le Cn |N_F(v)|^{\alpha-1}.
	\]
	Further, we can partition $N_R(v)$ into $r=\left\lceil \frac{|N_R(v)|}{|N_F(v)|}\right\rceil$ sets $X_1, \ldots, X_r$ with $|X_i|\le |N_F(v)|$ for every $i\in [r]$. Then,
	\[
	\begin{aligned}
	e(N_F(v), N_R(v)) = \sum_{i=1}^r e(N_F(v), X_i) &\le \sum_{i=1}^r \ex(2|N_F(v)|, H)\\& \le C\cdot r\cdot (2|N_F(v)|)^{\alpha}
	\\& \le C\cdot \frac{2|N_R(v)|}{|N_F(v)|} \cdot 2^{\alpha}|N_F(v)|^{\alpha}
	\\& \le 2^{1+\alpha} Cn|N_F(v)|^{\alpha-1}.
	\end{aligned}
	\]
	Therefore, (\ref{eq:concl1}) gives us
	\[
	t(G)\le \sum_{v\in V(G)}\left((2^{1+\alpha}+1)C \cdot n|N_F(v)|^{\alpha-1}\right) = C'n\sum_{v\in V(G)}|N_F(v)|^{\alpha-1},
	\]
	where $C' = (2^{1+\alpha}+1)C$. Observe that $\sum_{v\in V(G)}N_F(v) = 2|F| = 6\epsilon n^2$. Now, as $1<\alpha<2$, the function $f(x)=x^{\alpha -1}$ is concave in $x$, and thus by Jensen's inequality,
	\[
	t(G)\le C'n^2\cdot \left(\frac 1n \sum_{v\in V(G)}|N_F(v)|\right)^{\alpha-1} = C'n^2\cdot \left(6\epsilon n\right)^{\alpha-1} = C'' \epsilon^{\alpha-1} n^{1+\alpha},
	\]
	where $C'' = C'\cdot 6^{\alpha-1}$.
	As $\epsilon>0$ was arbitrary, this implies $t(G)=o(n^{1+\alpha})$.

\end{proof}
	
\bigskip

\paragraph{Acknowledgments.}
The first author is partially supported by NSF awards DMS-1300138, 1763317, 1952767 and 2153576. The second author is very thankful to Xizhi Liu for several helpful discussions and comments, and for sketching a proof of Lemma \ref{lem:pathsLemma}. We are also grateful to Jozsef Balogh, Cory Palmer and Abhishek Methuku for informing us about previous work on these problems. We would like to thank the anonymous referees for their helpful comments.

\bibliographystyle{plain}
\bibliography{suspensionFree.bbl}

\pagebreak
\section{Appendix}
Our goal in this section is to complete the proof of Lemma \ref{lem:cleaningForP5}. Recall that $H$ is a subgraph of a $\{K_6^-,\widehat P_5\}$-free graph $G$ such that every edge of $H$ has codegree at least $2$.
\begin{proof}[Proof of Lemma~\ref{lem:cleaningForP5}]
We wish to show that $H$ does not contain copies of $W_5$ or $K_{1,2,2}$. We do this via sequentially cleaning the following graphs from $H$:
\begin{itemize}
	\item $K_6^{-2,1}$, the graph obtained from $K_6$ by deleting two intersecting edges,
	\item $K_6^{-2,2}$, the graph obtained from $K_6$ by deleting two non-intersecting edges,
	\item $K_6^{-3,1}$, the graph obtained from $K_6$ by deleting a $P_3$,
	\item $K_6^{-3,2}$, the graph obtained from $K_6$ by deleting a $P_2\sqcup K_2$,
	\item $K_5$,
	\item $K_5^-$, the graph obtained from $K_5$ by deleting one edge,
	\item $W_5^+$, the graph obtained from the $5$-wheel graph $W_5=\widehat{C}_5$ by adding an edge,
	\item $W_5$, and
	\item $K_{1,2,2}$, the $4$-wheel graph.
\end{itemize}
More specifically, whenever $H$ contains a copy of one of these graphs, we would be able to use the induction hypothesis of (\ref{eq:inductionP5}) on some subgraph $H'\subsetneq H$ and complete the induction step for $H$.

Before proceeding onto the cleaning steps, we make an important observation:
\begin{equation}
\label{eq:appendixEqP4extend}
\mbox{If }abcde\mbox{ is a }P_4\mbox{ in }N_H(x),\mbox{ then } N_H(\{a,x\}) \subseteq \{b,c,d,e\}.
\end{equation}

This is because since $xa$ has codegree at least $2$, we must have a vertex $y\in V(H)$ with $xay$ being a triangle, $y\neq b$. If $y\not\in\{c,d,e\}$, then we get a $\widehat P_5$ centered around $x$. Thus $y=c$ or $y=d$ or $y=e$, implying (\ref{eq:appendixEqP4extend}).

\paragraph{1. Cleaning $K_6^{-2,1}$:}
Suppose $H$ has a copy of $K_6^{-2,1}$ with vertices $\{a,b,c,d,e,f\}$ such that the edges $ab$ and $bc$ are missing. This is an induced subgraph as $H$ has no $K_6^-$. Further, all edges of this subgraph other than $ac$ cannot belong to external triangles, as verified in Figure~\ref{fig:appendixFig1}.
\begin{figure}[h]
	\centering
	\begin{tikzpicture}
	\tikzstyle{vertex}=[circle,fill=black,minimum size=2pt,inner sep=1.3pt]
	\node[vertex](a1)at(1,0){};
	\node[vertex](a2)at(60:1){};
	\node[vertex](a3)at(120:1){};
	\node[vertex](a4)at(180:1){};
	\node[vertex](a5)at(240:1){};
	\node[vertex](a6)at(300:1){};
	\node[vertex](x)at(150:1.732){};
	\draw (a3)--(x)--(a4);
	\draw(a1)--(a3)--(a4)--(a5)--(a6)--(a2)--(a5)--(a1)--(a4)--(a2)--(a5);
	\draw(a1)--(a6)--(a3)--(a5);
	\draw(a6)--(a4);
	\draw(a4)circle(0.2);
	\draw [line width=1pt](x)--(a3)--(a1)--(a6)--(a5)--(a2);
	
	\begin{scope}[shift={(4,0)}]
	\node[vertex](a1)at(1,0){};
	\node[vertex](a2)at(60:1){};
	\node[vertex](a3)at(120:1){};
	\node[vertex](a4)at(180:1){};
	\node[vertex](a5)at(240:1){};
	\node[vertex](a6)at(300:1){};
	\node[vertex](x)at(210:1.732){};
	\draw (a5)--(x)--(a4);
	\draw(a1)--(a3)--(a4)--(a5)--(a6)--(a2)--(a5)--(a1)--(a4)--(a2)--(a5);
	\draw(a1)--(a6)--(a3)--(a5);
	\draw(a6)--(a4);
	\draw(a4)circle(0.2);
	\draw [line width=1pt](x)--(a5)--(a3)--(a1)--(a6)--(a2);
	\end{scope}
	\end{tikzpicture}
	\caption{All edges but $ac$ cannot lie in external triangles.\label{fig:appendixFig1}}
\end{figure}
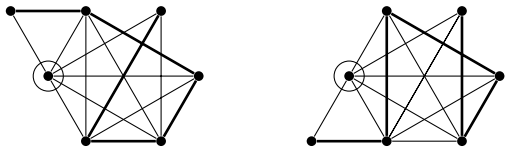

Now suppose $ac$ lies on an external triangle, $acx$. By (\ref{eq:appendixEqP4extend}) on the $4$-path $fedcx$ in $N_H(a)$, $N_H(\{a,x\})\subseteq\{c,d,e\}$. Moreover, $ax$ has codegree at least $2$. Thus, $xd$, $xe$, or $xf$ is an edge of $H$. In either case we obtain $\widehat P_5$'s in $H$, as shown in Figure~\ref{fig:appendixFig2}.
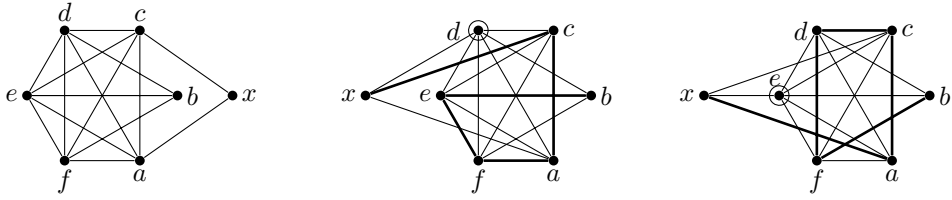
\begin{figure}[h]
	\centering
	\begin{tikzpicture}
	\tikzstyle{vertex}=[circle,fill=black,minimum size=2pt,inner sep=1.3pt]
	\begin{scope}[rotate=-60]
	\node[vertex](a1)at(1,0){};
	\node[vertex](a2)at(60:1){};
	\node[vertex](a3)at(120:1){};
	\node[vertex](a4)at(180:1){};
	\node[vertex](a5)at(240:1){};
	\node[vertex](a6)at(300:1){};
	\node[vertex](x)at(60:1.732){};
	\draw(a1)node[below]{$a$};
	\draw(a2)node[right]{$b$};
	\draw(a3)node[above]{$c$};
	\draw(a4)node[above]{$d$};
	\draw(a5)node[left]{$e$};
	\draw(a6)node[below]{$f$};
	\draw(x)node[right]{$x$};
	\draw (a3)--(x)--(a1);
	\draw(a1)--(a3)--(a4)--(a5)--(a6)--(a2)--(a5)--(a1)--(a4)--(a2)--(a5);
	\draw(a1)--(a6)--(a3)--(a5);
	\draw(a6)--(a4);
	\end{scope}
	
	\begin{scope}[shift={(5.5,0)},rotate=-60]
	\node[vertex](a1)at(1,0){};
	\node[vertex](a2)at(60:1){};
	\node[vertex](a3)at(120:1){};
	\node[vertex](a4)at(180:1){};
	\node[vertex](a5)at(240:1){};
	\node[vertex](a6)at(300:1){};
	\node[vertex](x)at(60:-2){};
	\draw(a1)node[below]{$a$};
	\draw(a2)node[right]{$b$};
	\draw(a3)node[right]{$c$};
	\draw(a4)node[left]{$d\ $};
	\draw(a5)node[left]{$e$};
	\draw(a6)node[below]{$f$};
	\draw(x)node[left]{$x$};
	\draw (a3)--(x)--(a1);
	\draw(a1)--(a3)--(a4)--(a5)--(a6)--(a2)--(a5)--(a1)--(a4)--(a2)--(a5);
	\draw(a1)--(a6)--(a3)--(a5);
	\draw(a6)--(a4);
	\draw(x)--(a4);
	\draw(a4)circle(0.13);
	\draw[line width=1pt](x)--(a3)--(a1)--(a6)--(a5)--(a2);
	\end{scope}
	
	\begin{scope}[shift={(10,0)},rotate=-60]
	\node[vertex](a1)at(1,0){};
	\node[vertex](a2)at(60:1){};
	\node[vertex](a3)at(120:1){};
	\node[vertex](a4)at(180:1){};
	\node[vertex](a5)at(240:1){};
	\node[vertex](a6)at(300:1){};
	\node[vertex](x)at(60:-2){};
	\draw(a1)node[below]{$a$};
	\draw(a2)node[right]{$b$};
	\draw(a3)node[right]{$c$};
	\draw(a4)node[left]{$d$};
	\draw(a5)node[above]{$e\ $};
	\draw(a6)node[below]{$f$};
	\draw(x)node[left]{$x$};
	\draw (a3)--(x)--(a1);
	\draw(a1)--(a3)--(a4)--(a5)--(a6)--(a2)--(a5)--(a1)--(a4)--(a2)--(a5);
	\draw(a1)--(a6)--(a3)--(a5);
	\draw(a6)--(a4);
	\draw(x)--(a5);
	\draw(a5)circle(0.13);
	\draw[line width=1pt](x)--(a1)--(a3)--(a4)--(a6)--(a2);
	\end{scope}
	\end{tikzpicture}
	\caption{$acx$ is a triangle.\label{fig:appendixFig2}}
\end{figure}

Hence $K_6^{-2,1}$ is a block by itself. Let $H'$ be the subgraph of $H$ obtained by deleting all edges from this copy of $K_6^{-2,1}$. Then note that $t(H')=t(H)-13$ and $e(H')=e(H)-13$, completing the induction step for $H$.

\paragraph{2. Cleaning $K_6^{-2,2}$:}
Let $H$ have a copy of $K_6^{-2,2}$ with vertices $\{a,b,c,d,e,f\}$ such that the edges $ab$ and $cd$ are missing. Clearly this is an induced subgraph of $H$. It can be checked that the edges $ea,ec,eb,ed$ cannot lie on external triangles as otherwise we would get $\widehat P_5$'s centered at $e$. Similarly, the edges $fc,fa,fd,fb$ cannot lie on external triangles. 

Now, suppose the edge $ad$ lies in an external triangle $adx$. Refer to Figure~\ref{fig:appendixFig3}. As $cfedx$ is a $P_4$ in the neighborhood of $a$, (\ref{eq:appendixEqP4extend}) implies that $N_H(\{a,x\})\subseteq \{d,e,f,c\}$ and $N_H(\{a,c\})\subseteq \{f,e,d,x\}$. If $xe\in E(H)$, this leads to a $\widehat P_5$ centered at $e$, given by the $5$-path $caxdbf$. If $xf\in E(H)$, we get a $\widehat P_5$ centered at $f$, given by the $5$-path $bdxac$. Thus, $xc\in E(H)$, and the edge $xa$ is light. Repeating the same argument for the $4$-path $bfeax$ around $d$, we get $xb\in E(H)$, $xd$ is light and $ac$ has codegree $3$. Using (\ref{eq:appendixEqP4extend}) on the $4$-path $xaefb$ in $N_H(c)$ and $N_H(d)$ respectively, we get $N_H(\{c,x\})\subseteq\{a,e,f,b\}$ and $N_H(\{b,d\})\subseteq\{f,e,a,x\}$. Since we already know that $xe,xf,ab\not\in E(H)$, this means that the edge $xc$ is light and $bd$ has codegree $3$. Similarly, $bx$ is light. Now, let
\[
H'=H-\{ea,ec,eb,ed,fc,fa,fd,fb,xa,xc,xb,xd,ac,bd\}.
\]
Clearly $e(H')=e(H)-14$ and $t(H')=t(H)-14$, and we are done by induction.

\begin{figure}[h]
	\centering
	\begin{tikzpicture}[scale=1.5]
	\tikzstyle{vertex}=[circle,fill=black,minimum size=2pt,inner sep=1.3pt]
	\node[vertex](a)at(0,1){};
	\node[vertex](c)at(0,0){};
	\node[vertex](e)at(1,1){};
	\node[vertex](f)at(1,0){};
	\node[vertex](d)at(2,1){};
	\node[vertex](b)at(2,0){};
	\node[vertex](x)at(1,2){};
	\draw (a)node[left]{$a$};
	\draw (b)node[right]{$b$};
	\draw (c)node[left]{$c$};
	\draw (d)node[right]{$d$};
	\draw (e)node[above]{$e$};
	\draw (f)node[below]{$f$};
	\draw (x)node[above]{$x$};
	\draw (b)--(c)--(a)--(d)--(b)--(e)--(c);
	\draw (a)--(f)--(d);
	\draw (e)--(f);
	\draw (a) to [bend left](d);
	\draw (c) to [bend right](b);
	\draw (a)--(x)--(d);
	\draw [dashed](c)--(x)--(b);
	\end{tikzpicture}
	\caption{Edge $ad$ lies in triangle $adx$, $x\not\in\{a,b,c,d,e,f\}$.\label{fig:appendixFig3}}
\end{figure}
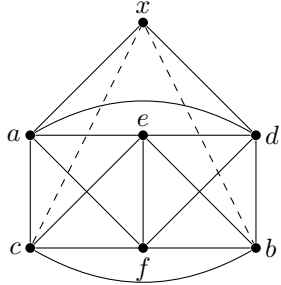

Therefore, $ad$ cannot lie in any external triangle $adx$, and is light. Similarly, $bc$ is light. Let $H'=H-\{ea,ec,eb,ed,fc,fa,fd,fb,ad,bc\}$. Then $t(H')=t(H)-10$ and $e(H')=e(H)-10$, finishing the induction step for $H$.

\paragraph{3. Cleaning $K_6^{-3,1}$:} Since $G$ has no $K_6^{-2,1}$ or $K_6^{-2,2}$ which are the only two ways one can delete two edges from $K_6$, any copy of $K_6^{-3,1}$ is induced. Suppose such a copy of $K_6^{-3,1}$ exists in $G$, and is given by the complete graph on $\{a,b,c,d,e,f\}$ minus the edges $\{ab,bc,cd\}$. By an argument exactly the same as before, $ea,ec,eb,ed,fc,fa,fd,fb$ are light. Further, if $ad$ lies in an external triangle $adx$ (as in Figure~\ref{fig:appendixFig4}), then by repeating the argument for cleaning $K_6^{-2,2}$, we note that $xc,xb\in E(H)$, $xa,xd$ are light, and $ac,bd$ have codegree exactly three. Also, by using (\ref{eq:appendixEqP4extend}) on the $4$-path $dxcfe$, we have $N_H(\{a,d\})\subseteq\{x,c,f,e\}$. As $cd\not\in E(H)$, this means $N_H(\{a,d\})=\{x,f,e\}$, and therefore $ad$ has codegree three. Finally, (\ref{eq:appendixEqP4extend}) on the path $cfedx$ in $N_H(a)$ gives us $N_H(\{a,c\})\subseteq \{f,e,d,x\}$, whereas $cd\not\in E(H)$, implying that $ac$ has codegree three as well. Similarly, $bd$ has codegree three.

Let $H'=H-\{ea,ec,eb,ed,fc,fa,fd,fb,ac,bd,ad,xa,xd\}$. Then, $t(H')=t(H)-13$ and $e(H')=e(H)-13$, and we can proceed by the induction hypothesis on $H'$.

\begin{figure}[h]
	\centering
	\begin{tikzpicture}[scale=1.5]
	\tikzstyle{vertex}=[circle,fill=black,minimum size=2pt,inner sep=1.3pt]
	\node[vertex](a)at(0,1){};
	\node[vertex](c)at(0,0){};
	\node[vertex](e)at(1,1){};
	\node[vertex](f)at(1,0){};
	\node[vertex](d)at(2,1){};
	\node[vertex](b)at(2,0){};
	\node[vertex](x)at(1,2){};
	\draw (a)node[left]{$a$};
	\draw (b)node[right]{$b$};
	\draw (c)node[left]{$c$};
	\draw (d)node[right]{$d$};
	\draw (e)node[above]{$e$};
	\draw (f)node[below]{$f$};
	\draw (x)node[above]{$x$};
	\draw (b)--(c)--(a)--(d)--(b)--(e)--(c);
	\draw (a)--(f)--(d);
	\draw (e)--(f);
	\draw (a) to [bend left](d);
	\draw (a)--(x)--(d);
	\draw [dashed](c)--(x)--(b);
	\end{tikzpicture}
	\caption{$adx$ is an external triangle, $ab,bc,ca$ are non-edges in $H$. \label{fig:appendixFig4}}
\end{figure}
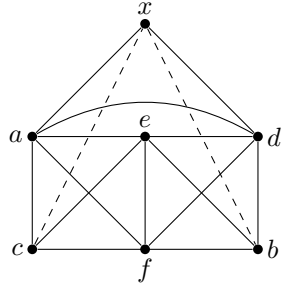

On the other hand, if the edge $ad$ is light, then we can simply let $H'=H-\{ea,ec,eb,ed,fc,fa,fd,fb,ad\}$, whence $t(H')=t(H)-9$ and $e(H')=e(H)-9$, and the induction step would be complete. Hence we can assume that $H$ is $K_6^{-3,1}$-free.

\paragraph{4. Cleaning $K_6^{-3,2}$:}
Suppose $H$ contains a $K_6^{-3,2}$ on vertices $\{a,b,c,d,e,f\}$ such that edges $ab,cd,de$ are missing. Since $H$ is $K_6^{-2,1}$ and $K_6^{-2,2}$-free, this subgraph is induced. As the edges $bd$ and $ad$ must have codegree at least two, there exist $x,y\in V(H)\setminus\{a,b,c,d,e,f\}$ such that $bdx$ and $ady$ are triangles in $H$. We consider two different cases.
\begin{itemize}
	\item {\bf Case 1. $x=y$} (Figure~\ref{fig:appendixFig4.5} (left)):
	Since $N_H(b)$ contains the 4-path $cefdx$, (\ref{eq:appendixEqP4extend}) gives $N_H(\{b,x\})\subseteq \{d,f,e,c\}$. If $xf\in E(H)$, then $N_H(f)$ contains the 5-path $xadbce$. If both $xc$ and $xe$ were edges in $H$, then $H[\{a,b,c,e,f,x\}]$ would be a $K_6^{-2,2}$ with edges $xf,ab$ missing. Therefore, only one of $xc$ and $xe$ can be an edge. By symmetry, assume $xc\in E(H)$ and $xe\not\in E(H)$.
	
	As this fixes edges and non-edges between any pair of vertices from $\{a,b,c,d,e,f,x\}$, $H[\{a,b,c,d,e,f,x\}]$ is induced. Consider the $5$-wheel $(f,adbcea)$, where the first tuple denotes the central vertex and the second tuple is the outer $C_5$. Since none of the edges $fa,fb,fc,fd,fe$ can lie in triangles with a  vertex $y\not\in \{a,b,c,d,e,f,x\}$ (it would give a $\widehat P_5$ around $f$), they all have exhausted their codegrees. Similarly, $(c,befaxb)$, $(b,cxdfec)$, and $(a,dxcefd)$ are $W_5$'s in $H$. Let
	\[
	H'=H-\{cb,cx,ca,cf,ce,fe,fb,fd,fa,be,bx,bd,ad,ax,ae\}.
	\]
	Then, $e(H')=e(H)-15$ and $t(H')=t(H)-13$ ($4$ triangles through $x$, $7$ through $f$ but not $x$, and $2$ not through $x$ or $f$). We can then proceed with the induction hypothesis on $H'$.
	\begin{figure}[h]
		\centering
		\begin{tikzpicture}[scale=1.5]
		\tikzstyle{vertex}=[circle,fill=black,minimum size=2pt,inner sep=1.3pt]
		\node[vertex](a)at(-30:1){};
		\node[vertex](b)at(100:1){};
		\node[vertex](c)at(150:1){};
		\node[vertex](d)at(30:1){};
		\node[vertex](e)at(210:1){};
		\node[vertex](f)at(260:1){};
		\node[vertex](x)at(30:0.3){};
		\draw (a)node[right]{$a$};
		\draw (b)node[above]{$b$};
		\draw (c)node[left]{$c$};
		\draw (d)node[right]{$d$};
		\draw (e)node[left]{$e$};
		\draw (f)node[below]{$f$};
		\draw (x)node[below]{$x\ $};
		\draw (c)--(a)--(d)--(f)--(c)--(b)--(e)--(c)--(f)--(a)--(e)--(f)--(b)--(d);
		\draw (a)--(x)--(d);
		\draw (x)--(b);
		\draw [dashed](x)--(c);
		
		\begin{scope}[shift={(4,0)}]
		\node[vertex](a)at(-30:1){};
		\node[vertex](b)at(100:1){};
		\node[vertex](c)at(150:1){};
		\node[vertex](d)at(30:1){};
		\node[vertex](e)at(210:1){};
		\node[vertex](f)at(260:1){};
		\node[vertex](x)at(65:0.3){};
		\node[vertex](y)at(0:1.4){};
		\draw (a)node[below]{$a$};
		\draw (b)node[above]{$b$};
		\draw (c)node[left]{$c$};
		\draw (d)node[above]{$d$};
		\draw (e)node[left]{$e$};
		\draw (f)node[below]{$f$};
		\draw (x)node[below]{$x\ $};
		\draw (y)node[right]{$y$};
		\draw (c)--(a)--(d)--(f)--(c)--(b)--(e)--(c)--(f)--(a)--(e)--(f)--(b)--(d);
		\draw (b)--(x)--(d);
		\draw (a)--(y)--(d);
		\draw [dashed](x)--(c);
		\draw [dashed](y)--(e);
		\end{scope}
		\end{tikzpicture}
		\caption{$bdx$ and $ady$ are triangles, (left: $x=y$, right: $x\ne y$).\label{fig:appendixFig4.5}}
	\end{figure}
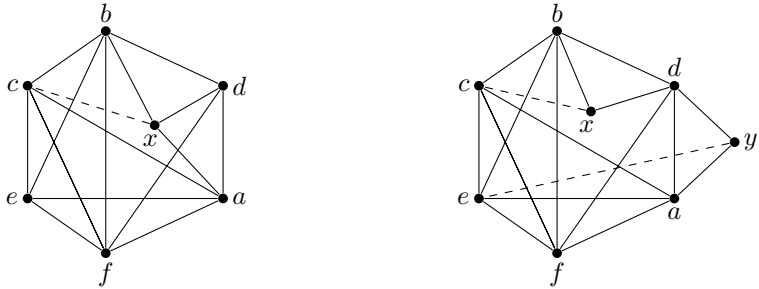

	\item {\bf Case 2. $x\neq y$} (Figure~\ref{fig:appendixFig4.5} (right)):
	Without loss of generality assume $by,ax\not\in E(H)$, as these would lead us to Case 1. As $N_H(b)$ contains the $4$-path $xdfec$, by (\ref{eq:appendixEqP4extend}), $N_H(\{b,x\})\subseteq\{d,f,e,c\}$. Note that if $xf\in E(H)$, then $N_H(f)$ contains the path $ecbxda$ of length $5$. Hence, $N_H(\{b,x\})\subseteq\{d,e,c\}$. Further, both $xc$ and $xe$ cannot be edges in $H$, as then $H[\{a,b,c,e,f,x\}]\supseteq K_6^{-3,1}$ with the edges $ba,ax,xf$ missing. As codegree of $bx$ is at least $2$, exactly one of $xc$ and $xe$ is an edge in $H$. By symmetry, assume $xc\in E(H)$ and $xe\not\in E(H)$.

	Now by (\ref{eq:appendixEqP4extend}) on the $4$-path $ydfec$ in $N_H(a)$, we get $N_H(\{a,y\})\subseteq \{d,f,e,c\}$. If $yf\in E(H)$, then $N_H(f)$ contains the path $aydbce$ of length $5$, and if $yc\in E(H)$, then $N_H(c)$ contains the path $xbefay$ of length $5$. Therefore, $N_H(\{a,y\})=\{d,e\}$, and $ye\in E(H)$.
	
	Finally, let us consider the $4$-path $yafbc$ in $N_H(e)$. Using (\ref{eq:appendixEqP4extend}),
	\[
	N_H(\{y,e\})\subseteq \{a,f,b,c\}.
	\]
	However, $yf,yc\not\in E(H)$ from our argument in the last paragraph, and $by\not\in E(H)$ as we are in Case 2. This is a contradiction, as the edge $ye$ must have codegree at least $2$.
\end{itemize}

\paragraph{5. Cleaning $K_5$:}
If $H$ contains a copy of $K_5$ on vertex set $\{a,b,c,d,e\}$, then we claim that it is a block by itself. Suppose $x\in V(H)\setminus\{a,b,c,d,e\}$ is such that $abx$ is a triangle in $H$. Since $xbcde$ is a $P_4$ in $N_H(a)$, (\ref{eq:appendixEqP4extend}) implies that $N_H(\{a,x\})\subseteq \{b,c,d,e\}$. Further, $ax$ must have codegree at least $2$. Thus, $xc$, $xd$, or $xe$ is an edge. In either case, $H[\{a,b,c,d,e,x\}]\supseteq K_6^{-2,1}$, a contradiction.

\paragraph{6. Cleaning $K_5^-$:}
Let $H$ have a copy of $K_5^-$ on vertices $a,b,c,d,e$ such that $ab\not\in E(H)$. If the edge $bc$ lies in an external triangle $bcx$ as shown in Figure~\ref{fig:appendixFig5}, then note that $xb$ has codegree at least two, and (\ref{eq:appendixEqP4extend}) on the $4$-path $xbdea$ in $N_H(c)$ tells us that $N_H(\{c,x\})\subseteq \{b,d,e,a\}$. If $xe\in E(H)$ then $G[a,b,c,d,e,x]$ contains the graph $K_6^{-3,1}$ with edges $dx,xa,ab$ missing. If $xd\in E(H)$, then we have the $K_6^{-3,1}$ with edges $ex,xa,ab$ missing. Finally, if $xa\in E(H)$, then $G[a,b,c,d,e,x]$ contains $K_6^{-3,2}$ with edges $ex,xd,ab$ missing. Thus, the edge $bc$ cannot lie on an external triangle.

\begin{figure}[h]
	\centering
	\begin{tikzpicture}
	\tikzstyle{vertex}=[circle,fill=black,minimum size=2pt,inner sep=1.3pt]
	\begin{scope}[rotate=234]
	\node[vertex](a)at(1,0){};
	\node[vertex](b)at(72:1){};
	\node[vertex](c)at(144:1){};
	\node[vertex](d)at(216:1){};
	\node[vertex](e)at(288:1){};
	\node[vertex](x)at(108:2){};
	\draw(b)--(c)--(d)--(e)--(a)--(c)--(e)--(b)--(d)--(a);
	\draw(b)--(x)--(c);
	\draw(a)node[below]{$a$};
	\draw(b)node[below]{$b$};
	\draw(c)node[above]{$\ \ c$};
	\draw(d)node[above]{$d$};
	\draw(e)node[left]{$e$};
	\draw(x)node[right]{$x$};
	\end{scope}

	\begin{scope}[shift={(5,0)},rotate=162]
	\node[vertex](a)at(1,0){};
	\node[vertex](b)at(72:1){};
	\node[vertex](c)at(144:1){};
	\node[vertex](d)at(216:1){};
	\node[vertex](e)at(288:1){};
	\node[vertex](x)at(216:2){};
	\draw(b)--(c)--(d)--(e)--(a)--(c)--(e)--(b)--(d)--(a);
	\draw(d)--(x)--(c);
	\draw(a)node[below]{$a\,$};
	\draw(b)node[below]{$b$};
	\draw(c)node[below]{$c$};
	\draw(d)node[below]{$\, d$};
	\draw(e)node[above]{$e$};
	\draw(x)node[right]{$x$};
	\draw (x)--(e);
	\end{scope}
	\end{tikzpicture}
	\caption{$K_5^-$ is a block by itself.\label{fig:appendixFig5}}
\end{figure}

Thus by symmetry, $ae,ad,ac,be,bd,bc$ cannot lie on external triangles. Now suppose that the edge $cd$ lies on an external triangle $cdx$. By (\ref{eq:appendixEqP4extend}) on any $\widehat P_4$ centered at $c$, $N_H(\{c,x\})\subseteq \{d,e,a,b\}$. If either $xa$ or $xb$ is an edge, we obtain a $K_6^{-3,1}$ with missing edges $ab,bx,xe$ or $ba,ax,xe$, respectively. So assume $xa,xb\not\in E(H)$. Thus $ex\in E(H)$, and $cx$ has codegree $2$. Similarly, $dx$ has codegree $2$. Now using (\ref{eq:appendixEqP4extend}) on the $4$-path $xdbca$ in $N_H(e)$, we have $N_H(\{e,x\})\subseteq \{a,b,c,d\}$. Since $xa,xb\not\in E(H)$, $ex$ must have codegree $2$. Thus, the edges $xc,xd,xe$ all have codegree $2$. Let $H'=H-\{xc,xd,xe\}$, then $t(H')=t(H)-3$ and $e(H')=e(H)-3$, and we can proceed by induction.

Hence, we may assume that $cd$ also does not lie on external triangles. Let $H'=H-\{ac,ad,ae,bc,bd,be,cd\}$. Then, $t(H')=t(H)-7$ and $e(H')=e(H)-7$, completing the induction hypothesis for $H$ again. So, without loss of generality we can assume that $G$ is $\{K_5^-,\widehat P_5\}$-free. 

\paragraph{7. Cleaning $W_5^+$:}
Let $H$ contain a $W_5^+$, given by central vertex $x$, outer cycle $abcde$ with an edge $ac\in E(H)$. If $ad\in E(H)$, then $a,b,c,d,x$ form a $K_5^-$ in $H$. Therefore by symmetry, all copies of $W_5^+$ in $H$ are induced.

Now let us fix such a $W_5^+$ in $H$ with vertices labeled as above. As $cd$ and $ae$ must have codegree at least $2$, there exist $y,z\in V(H)\setminus \{a,b,c,d,e,x\}$ such that $aey$ and $cdz$ are triangles in $H$. Then, we have two possible cases:
\begin{itemize}
	\item {\bf Case 1:} $y=z$. 
	Refer to Figure~\ref{fig:appendixFig6}. Note that if $yx\in E(H)$, then $N_H(x)$ would contain the $P_5$ given by $yabcde$. Hence $yx\not\in E(H)$. Using (\ref{eq:appendixEqP4extend}) on the $4$-path $yexcb$, $N_H(\{a,y\})\subseteq \{e,x,c,b\}$. Suppose $yb\in E(H)$. Note that $(y,abcdea)$, $(c,aydxba)$, $(a,cybxec)$ form $W_5^+$'s in $H$. As (\ref{eq:appendixEqP4extend}) together with the fact that every $W_5^+$ of $H$ is induced imply that each central edge of any copy of $W_5^+$ in $H$ does not lie on external triangles, the edges $xa,xb,xc,xd,xe$; $ya,yb,yc,yd,ye$; $cb,ca,cd$; $ab,ae$ cannot lie in external triangles. Thus, we can delete these $15$ edges from $H$, and only lose $13$ triangles ($6$ through $x$, $6$ through $y$, and the triangle $abc$). Our proof would then be complete by induction.
	
	Hence, assume $yb\not\in E(H)$. From $N_H(\{a,y\})\subseteq \{e,x,c,b\}$ this implies that codegree of $ay$ is exactly $2$. Similarly, $cy$ has codegree $2$. Further, (\ref{eq:appendixEqP4extend}) on the path of length four $bxcye$ in $N_H(a)$ implies that $N_H(\{a,b\})\subseteq \{x,c,y,e\}$. Since $by,be\not\in E(H)$, this implies that $N_H(\{a,b\})=\{x,c\}$, and $ab$ has codegree $2$. Similarly, $bc$ has codegree $2$. Let $H'=H-\{xa,xb,xc,xd,xe,ab,bc\}$. Then, $t(H')=t(H)-7$ and $e(H')=e(H)-7$, again concluding the induction step.

	\begin{figure}[h]
		\centering
		\begin{tikzpicture}
		\tikzstyle{vertex}=[circle,fill=black,minimum size=2pt,inner sep=1.3pt]
		\begin{scope}[rotate=108]
		\node[vertex](a)at(1,0){};
		\node[vertex](b)at(72:1){};
		\node[vertex](c)at(144:1){};
		\node[vertex](d)at(216:1){};
		\node[vertex](e)at(288:1){};
		\node[vertex](x)at(0,0){};
		\node[vertex](y)at(252:2){};
		\draw(a)--(b)--(c)--(d)--(e)--(a)--(x)--(c)--(x)--(d);
		\draw(e)--(x)--(b);
		\draw(a)--(c);
		\draw(e)--(y)--(d);
		\draw (a)to[bend left=45](y);
		\draw (y)to[bend left=45](c);
		\draw (y)to[bend right=20](b);
		\draw(a)node[above]{$a\ \ $};
		\draw(b)node[left]{$b$};
		\draw(c)node[below]{$c\ \ $};
		\draw(d)node[below]{$d$};
		\draw(e)node[above]{$e$};
		\draw(x)node[right]{$\ x$};
		\draw(y)node[right]{$y$};
		\end{scope}
		
		\begin{scope}[shift={(6,0)},rotate=108]
		\node[vertex](a)at(1,0){};
		\node[vertex](b)at(72:1){};
		\node[vertex](c)at(144:1){};
		\node[vertex](d)at(216:1){};
		\node[vertex](e)at(288:1){};
		\node[vertex](x)at(0,0){};
		\node[vertex](y)at(252:2){};
		\draw(a)--(b)--(c)--(d)--(e)--(a)--(x)--(c)--(x)--(d);
		\draw(e)--(x)--(b);
		\draw(a)--(c);
		\draw(e)--(y)--(d);
		\draw (a)to[bend left=45](y);
		\draw (y)to[bend left=45](c);
		\draw(a)node[above]{$a\ \ $};
		\draw(b)node[left]{$b$};
		\draw(c)node[below]{$c\ \ $};
		\draw(d)node[below]{$d$};
		\draw(e)node[above]{$e$};
		\draw(x)node[right]{$\ x$};
		\draw(y)node[right]{$y$};
		\end{scope}
		\end{tikzpicture}
		\caption{Case 1: $y=z$, $aey$ and $cdy$ are triangles in $H$.\label{fig:appendixFig6}}
	\end{figure}
	
	\item {\bf Case 2:} $y\neq z$.
	Since $N_H(a)$ contains the path $yexcb$ of length $4$, we must have $N_H(\{a,y\})\subseteq \{e,x,b,c\}$. If $yx\in E(H)$, we obtain the $5$-path $yabcde$ in $N_H(x)$, and if $yc\in E(H)$ we obtain the $5$-path $yabxdz$ in $N_H(c)$. On the other hand, $ay$ must have codegree at least $2$. Thus $N_H(\{a,y\})=\{b,e\}$ and $ay$ is light. By a symmetric argument, $zb\in E(H)$ and $cz$ is also light. Refer to Figure~\ref{fig:appendixFig7}. Since $N_H(b)$ contains the $P_4$ given by $zcxay$, we have $N_H(\{b,y\})\subseteq \{a,x,c,z\}$. However, we have already observed that $yx,yc\not\in E(H)$. Hence $zy\in E(H)$, and $by$ is light. Similarly, $bz$ is light.

	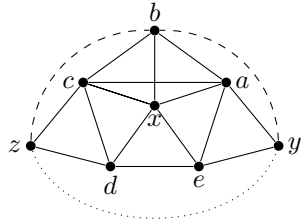
\begin{figure}[h]
		\centering
		\begin{tikzpicture}
		\tikzstyle{vertex}=[circle,fill=black,minimum size=2pt,inner sep=1.3pt]
		\begin{scope}[rotate=18]
		\node[vertex](a)at(1,0){};
		\node[vertex](b)at(72:1){};
		\node[vertex](c)at(144:1){};
		\node[vertex](d)at(216:1){};
		\node[vertex](e)at(288:1){};
		\node[vertex](x)at(0,0){};
		\node[vertex](y)at(324:1.732){};
		\node[vertex](z)at(180:1.732){};
		\draw(a)--(b)--(c)--(d)--(e)--(a)--(x)--(c)--(x)--(d);
		\draw(e)--(x)--(b);
		\draw(a)--(c);
		\draw (a)--(y)--(e);
		\draw (c)--(z)--(d);
		\draw[dashed] (y)to[bend right=45](b);
		\draw[dashed] (z)to[bend left=45](b);
		\draw[dotted] (z)to[bend right=70](y);
		\draw(a)node[right]{$a$};
		\draw(b)node[above]{$b$};
		\draw(c)node[left]{$c$};
		\draw(d)node[below]{$d$};
		\draw(e)node[below]{$e$};
		\draw(x)node[below]{$x$};
		\draw(y)node[right]{$y$};
		\draw(z)node[left]{$z$};
		\end{scope}
		\end{tikzpicture}
		\caption{Case 2: $y\neq z$, $aey$ and $cdz$ are triangles in $H$.\label{fig:appendixFig7}}
	\end{figure}
	
	Now, observe that we have produced two $W_5^+$'s given by $(c,zdxabz)$ and $(a,yexcby)$, both with the extra edge $ac$. By (\ref{eq:appendixEqP4extend}) in $N_H(c)$ and $N_H(b)$, all of the central edges cannot lie in external triangles. Let
	\[
	H' = H-\{cz,cd,cx,cb,ay,ae,ax,ab,ac,by,bz\}.
	\]
	It is clear that $e(H')=e(H)-11$, and that deleting these edges, we delete $6$ triangles through $c$ and $4$ triangles through $a$ that do not contain $c$, and the triangle $byz$ through $b$ which does not contain $a$ or $c$. Hence, $t(H')=t(H)-11$. This completes our induction step.
\end{itemize}

We may therefore assume that $H$ is $W_5^+$-free.

\paragraph{8. Cleaning $W_5$:}
Suppose $H$ has a copy of $W_5$ given by $(x,abcdea)$. As $H$ is $W_5^+$-free, $H[\{a,b,c,d,e,x\}]\cong W_5$. By (\ref{eq:appendixEqP4extend}) applied to $N_H(x)$, every central edge is light. Thus, we may let $H'=H-\{xa,xb,xc,xd,xe\}$, whence $t(H')=t(H)-5$ and $e(H')=e(H)-5$, allowing us to complete the induction step for $H$.

\paragraph{9. Cleaning $K_{1,2,2}$:}
Finally, let $H$ contain a $K_{1,2,2}$ with central vertex $x$ and outer cycle $abcd$. Since $H$ is $K_5^-$-free, $H[\{a,b,c,d,x\}]\cong K_{1,2,2}$. We claim that none of the edges $xa,xb,xc,xd$ lie on an external triangle.

For the sake of contradiction, assume $y\in V(H)\setminus\{a,b,c,d,x\}$ is such that $xay$ is a triangle in $H$ (Figure~\ref{fig:appendixFig8}). By (\ref{eq:appendixEqP4extend}) in $N_H(a)$, we have $N_H(\{a,y\})\subseteq \{x,d,c,b\}$. If $yd\in E(H)$, we obtain the $5$-wheel $(x,ydcbay)$, and if $yb\in E(H)$, we obtain the $5$-wheel $(x,yadcby)$. Since $|N_H(\{a,y\})|\ge 2$, we must have $N_H(\{a,y\})=\{x,c\}$. But then, $H[\{c,x,b,y,a\}]\cong K_5^-$, a contradiction.

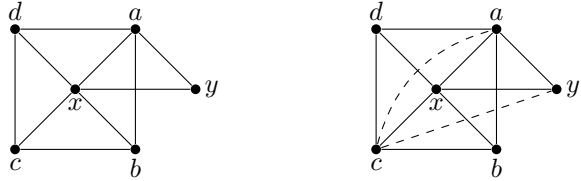
\begin{figure}[h]
	\centering
	\begin{tikzpicture}[scale=0.8]
	\tikzstyle{vertex}=[circle,fill=black,minimum size=2pt,inner sep=1.3pt]
	\node[vertex](a)at(1,1){};
	\node[vertex](b)at(1,-1){};
	\node[vertex](c)at(-1,-1){};
	\node[vertex](d)at(-1,1){};
	\node[vertex](x)at(0,0){};
	\node[vertex](y)at(2,0){};
	\draw(a)--(b)--(c)--(d)--(x)--(a)--(d);
	\draw(c)--(x)--(b);
	\draw (a)--(y)--(x);
	\draw(a)node[above]{$a$};
	\draw(b)node[below]{$b$};
	\draw(c)node[below]{$c$};
	\draw(d)node[above]{$d$};
	\draw(x)node[below]{$x$};
	\draw(y)node[right]{$y$};

	\begin{scope}[shift={(6,0)}]
	\node[vertex](a)at(1,1){};
	\node[vertex](b)at(1,-1){};
	\node[vertex](c)at(-1,-1){};
	\node[vertex](d)at(-1,1){};
	\node[vertex](x)at(0,0){};
	\node[vertex](y)at(2,0){};
	\draw(a)--(b)--(c)--(d)--(x)--(a)--(d);
	\draw(c)--(x)--(b);
	\draw (a)--(y)--(x);
	\draw[dashed] (c)--(y);
	\draw[dashed] (c)to[bend left](a);
	\draw(a)node[above]{$a$};
	\draw(b)node[below]{$b$};
	\draw(c)node[below]{$c$};
	\draw(d)node[above]{$d$};
	\draw(x)node[below]{$x$};
	\draw(y)node[right]{$y$};
	\end{scope}
	\end{tikzpicture}
	\caption{$xa$ cannot lie in an external triangle $xay$. \label{fig:appendixFig8}}
\end{figure}

Hence, the edges $xa,xb,xc,xd$ all have codegree $2$. Let $H'=H-\{xa,xb,xc,xd\}$, then $t(H')=t(H)-4$ and $e(H')=e(H)-4$, finishing the induction step in this case as well.

\medskip
Hence, after these cleaning procedures, we may assume that $H$ is a $\{W_5,K_{1,2,2}\}$-free subgraph of $G$ such that every edge of $H$ has codegree at least $2$. This concludes the proof of Lemma~\ref{lem:cleaningForP5}.
\end{proof}

\end{document}